\documentclass{amsart}

\usepackage{amssymb}
\usepackage{amsthm}
\usepackage[all]{xy}
\usepackage{color}
\usepackage{graphicx}

\author{Helmut Lenzing}
\address{Institut f\"{u}r Mathematik, Universit\"{a}t Paderborn, Warburger Str. 100, 33098-Paderborn, Germany}
\email{helmut@math.uni-paderborn.de}
\author{Jos\'{e}-Antonio de la Pe\~{n}a}
\address{Instituto de Matem\'{a}ticas, Universidad Nacional Aut\'{o}noma de M\'{e}xico, Circuito Exterior, M\'{e}xico 04510 DF, M\'{e}xico}
\email{jap@matem.unam.mx}

\newtheorem{theorem}{Theorem}[section]

\newtheorem{lemma}[theorem]{Lemma}

\newtheorem{proposition}[theorem]{Proposition}
\newenvironment{Example}{\begin{example}\em}{\end{example}}
\newtheorem{example}[theorem]{Example}
\newtheorem{examples}[theorem]{Examples}
\newenvironment{Examples}{\begin{examples}\em}{\end{examples}}
\newtheorem*{coro}{Corollary} 
\newtheorem{definition}[theorem]{Definition}
\newtheorem{remark}[theorem]{Remark}

 \newcommand{\cox}[1]{\Phi_{#1}}
 \newcommand{\coxpol}[1]{\chi_{#1}}
 \newcommand{\extcoxpol}[1]{\hat\chi_{#1}}
 \newcommand{\charpol}[1]{\kappa_{#1}}
 \newcommand{\stcharpol}[1]{\kappa^*_{#1}}

 \newcommand{\Der}[1]{\operatorname{D^\mathit{b}}#1}
 
 \newcommand{\mmod}[1]{\operatorname{mod}\textrm{-} #1}
 \newcommand{\proj}[1]{\operatorname{proj}\textrm{-}#1}
 \newcommand{\modgrad}[2]{\operatorname{mod}^{#1}\textrm{-}#2}
 \newcommand{\modgradnull}[2]{\operatorname{mod}^{#1}_0\textrm{-}#2}

 \newcommand{\projgrad}[2]{\operatorname{proj}^{#1}\textrm{-}#2}
 \newcommand{\Hom}[3]{\operatorname{Hom}_{#1}(#2,#3)}
 \newcommand{\Homgrad}[4]{\operatorname{Hom}^{#1}_{#2}(#3,#4)}
 \newcommand{\ulHom}[3]{\underline{\mathrm{Hom}}_{#1}(#2,#3)}
 \newcommand{\CMgrad}[2]{\operatorname{CM}^{#1}(#2)}
 \newcommand{\stabCMgrad}[2]{\operatorname{\underline{CM}}^{#1}(#2)}
 \newcommand{\dual}[1]{\operatorname{D}\,#1}
 \newcommand{\Knull}[1]{\operatorname{K}_0(#1)}
 \newcommand{\define}[1]{\emph{#1}}
 \newcommand{\transpose}[1]{{#1}^{t}}
 \newcommand{\ddim}[2]{{\operatorname{dim}_{#1}#2}}
 \newcommand{\homol}[2]{\operatorname{H}^{#1}(#2)}
 \newcommand{\ra}{\rightarrow}
 \newcommand{\eulerform}[2]{\langle #1,#2 \rangle}
 \newcommand{\eulerchar}[1]{\chi_{#1}}
 \newcommand{\ZZ}{\mathbb{Z}}
 \newcommand{\EE}{\mathbb{E}}
 \newcommand{\ZZplus}{{\mathbb{Z}_+}}
 \newcommand{\Ext}[4]{\operatorname{Ext}^{#1}_{#2}(#3,#4)}
 \newcommand{\inverse}[1]{#1^{-1}}
 \newcommand{\specrad}[1]{\rho_{#1}}
 \newcommand{\totient}{\varphi}
 \newcommand{\FF}{\mathbb{F}}
 \newcommand{\Aext}[1]{\tilde{\mathbb{A}}_{#1}}
\newcommand{\Dext}[1]{\tilde{\mathbb{D}}_{#1}}
\newcommand{\Eext}[1]{\tilde{\mathbb{E}}_{#1}}
\newcommand{\A}[1]{\mathbb{A}_{#1}}
\newcommand{\B}[1]{\mathbb{B}_{#1}}
\newcommand{\C}[1]{\mathbb{C}_{#1}}
\newcommand{\D}[1]{\mathbb{D}_{#1}}
\newcommand{\E}[1]{\mathbb{E}_{#1}}

\newcommand{\coh}[1]{\operatorname{coh}(#1)}
\newcommand{\Oo}{\mathcal{O}}
\newcommand{\vx}{\vec{x}}
\newcommand{\vy}{\vec{y}}
\newcommand{\vc}{\vec{c}}
\newcommand{\Hh}{\mathcal{H}}
\newcommand{\lra}{\longrightarrow}

\renewcommand{\AA}{\mathbb{A}}
\newcommand{\XX}{\mathbb{X}}

\newcommand{\XXbar}{{\bar\XX}}
\newcommand{\RR}{\mathbb{R}}
\renewcommand{\SS}{\mathbb{S}}
\newcommand{\Roots}[1]{\operatorname{Roots}(#1)}
\newcommand{\Spec}[1]{\operatorname{Spec}(#1)}
\newcommand{\hoch}[2]{\operatorname{H}^{#1}#2}
\newcommand{\tief}[2]{\operatorname{H}_{#1}#2}
\newcommand{\gldim}[1]{\mathrm{gl}\textrm{.}\mathrm{dim}(#1)}
\newcommand{\La}{\Lambda}
\newcommand{\la}{\lambda}
\newcommand{\pP}{{\emph{$\boldsymbol{p}$}}}
\newcommand{\Lp}{\operatorname{L}(\pP)}
\newcommand{\End}[2]{\operatorname{End}_{#1}(#2)}
\newcommand{\Dsing}[1]{\operatorname{D}_{\textrm{Sg}}(#1)}

\newcommand{\Lahat}{\hat\La}
\newcommand{\LL}{\mathrm{L}}
\newcommand{\Dsinggrad}[2]{\operatorname{D}_{\textrm{Sg}}^{#1}(#2)}
\newcommand{\GL}[2]{\mathrm{GL}_{#1}(#2)}
\newcommand{\ci}{\circ}

\newcommand{\lala}{{\emph{$\boldsymbol{\la}$}}}
\newcommand{\up}[1]{\stackrel{#1}{\lra}}
\newcommand{\vDe}{\vec{\Delta}}
\newcommand{\vGa}{\vec{\Gamma}}
\newcommand{\Tt}{\mathcal{T}}
\newcommand{\dis}{\displaystyle}
\newcommand{\de}{\delta}
\newcommand{\DD}{\mathbb{D}}
\newcommand{\pla}{\pP,\lala}
\newcommand{\Lapla}{{\La(\pla)}}
\newcommand{\Spla}{\operatorname{S}(\pla)}
\newcommand{\wtilde}{\widetilde}
\newcommand{\Xpla}{\mathbb{X}(\pla)}
\newcommand{\vect}[1]{\operatorname{vect}#1}
\newcommand{\Rpla}{\operatorname{R}(\pla)}
\newcommand{\Dir}{\bigoplus}
\newcommand{\vom}{\vec{\omega}}
\newcommand{\Ll}{\mathcal{L}}

\newcommand{\spitz}[1]{\langle#1\rangle}

\newcommand{\Mm}{\mathcal{M}}
\newcommand{\Mmplus}{\mathcal{M}_+}
\newcommand{\Ppplus}{\mathcal{P}_+}
\newcommand{\Ssplus}{\mathcal{S}_+}
\newcommand{\Dd}{\mathcal{D}}
\newcommand{\Gaplus}{\Gamma_{+}}
\newcommand{\lperp}[1]{{{}^{\perp}#1}}
\newcommand{\rperp}[1]{{#1^{\perp}}}
\newcommand{\bR}{\operatorname{\mathbf{R}}}
\newcommand{\incl}{\hookrightarrow}
\newcommand{\Ga}{\Gamma}
\newcommand{\stabvect}[1]{\operatorname{\underline{\mathrm{vect}}}#1}
\newcommand{\iso}{\cong}
\newcommand{\ul}{\underline}
\newcommand{\Cc}{\mathcal{C}}
\newcommand{\De}{\Delta}
\newcommand{\QQ}{\mathbb{Q}}
\newcommand{\union}{\cup}
\newcommand{\Proj}[2]{\mathbb{P}_{#1}(#2)}
\newcommand{\set}[1]{\{#1\}}
\newcommand{\ddir}{\oplus}
\newcommand{\CC}{\mathbb{C}}
\newcommand{\ga}{\gamma}
\newcommand{\al}{\alpha}
\newcommand{\Modgrad}[2]{\operatorname{Mod}^{#1}\textrm{-}#2}
\newcommand{\be}{\beta}
\newcommand{\Ee}{\mathcal{E}}
\newcommand{\Ff}{\mathcal{F}}
\newcommand{\Ii}{\mathcal{I}}
\newcommand{\om}{\omega}
\newcommand{\ka}{\kappa}
\newcommand{\ofrac}[2]{\genfrac{}{}{0pt}{}{#1}{#2}}

\begin{document}

 \title[Spectral analysis and singularities]{Spectral analysis of finite dimensional \linebreak algebras and singularities}

\begin{abstract}
For a finite dimensional algebra $A$ of finite global dimension the bounded derived category of finite dimensional $A$-modules admits Auslander-Reiten triangles such that the Auslander-Reiten translation $\tau$ is an equivalence. On the level of the Grothendieck group $\tau$ induces the Coxeter transformation $\cox{A}$. More generally this extends to a homologically finite triangulated category $\Tt$ admitting Serre duality. In both cases the Coxeter polynomial, that is, the characteristic polynomial of the Coxeter transformation yields an important homological invariant of $A$ or $\Tt$. Spectral analysis is the study of this interplay, it often reveals unexpected links between apparently different subjects.

This paper gives a summary on spectral techniques and studies the links to singularity theory. In particular, it offers a contribution to the categorification of the Milnor lattice through triangulated categories which are naturally attached to a weighted projective line.
\end{abstract}

\subjclass{Primary 16G20, 16G50; Secondary 14G60, 18E30}

\keywords{Coxeter transformation,
hereditary algebra, canonical algebra, extended canonical algebra, supercanonical algebra, weighted projective line, triangulated category of singularities, stable category of vector bundles, Milnor lattice}

\maketitle

\section*{Introduction}
Let $k$ be a field and $A$ be a finite dimensional $k$-algebra. If $A$ has finite global dimension, a theorem of Happel~\cite{Happel:1988} asserts that the bounded derived category $\Der{A}=\Der{\mmod{A}}$ of finite dimensional $A$-modules has Serre duality in the form
$$
\dual{\Hom{}{X}{Y[1]}}=\Hom{}{Y}{\tau X}
$$
where $\tau$ is a self-equivalence of $\Der{A}$. In particular, $\Der{A}$ has almost-split triangles and the equivalence $\tau$ serves as the Auslander-Reiten translation of $\Der{A}$.

In the above setting the \define{Grothendieck group} $\Knull{A}$ of the category $\mmod{A}$ of finite dimensional (right) $A$-modules, formed with respect to short exact sequences, naturally agrees with the Grothendieck group $\Knull{\Der{A}}$ of the derived category, formed with respect to exact triangles. Hence the Auslander-Reiten translation of $\Der{A}$ induces an automorphism of $\Knull{A}$, called the \define{Coxeter transformation} $\cox{A}$ of $A$. Since $\Knull{A}$ is the free abelian group on the (isomorphism classes of) simple $A$-modules, we finally obtain the characteristic polynomial $\coxpol{A}$ of $\cox{A}$, called the \define{Coxeter polynomial} of $A$.

By construction, the Coxeter polynomial is a derived invariant of $A$, that is, two algebras $A$ and $B$ (of finite global dimension) with triangle-equivalent bounded derived categories have identical Coxeter polynomials. The set of roots of the Coxeter transformation equals the \define{spectrum of} the Coxeter transformation. The \define{spectral radius} of $\cox{A}$, that is, the maximum of the absolute values of the complex roots of $\coxpol{A}$ will also be called the spectral radius of $A$ or $\coxpol{A}$.

The main task of the present article is to investigate the impact of the spectrum of $A$ on the ring- and module-theoretic properties of $A$ and to identify important classes of finite dimensional algebras by their spectral properties. While it is true that algebras with very different features can have the same spectrum, under reasonable restrictions the information on $A$ provided by its spectrum is surprisingly good.

Further, the spectral analysis of finite dimensional algebras yields an important window linking representation theory of finite dimensional algebras to many different areas of current research.
The link is provided through the following observation. Let $P_1,\ldots,P_n$ be a complete, representative system of indecomposable projective $A$-modules. The integers
$
c_{ij}=\ddim{k}{\Hom{}{P_i}{P_j}}
$
yield an integral $n\times n$-matrix $C=(c_{ij})$ with the property that
$
\Phi=-C^{-1}C^t
$
represents the Coxeter transformation $\cox{A}$ in the basis of $\Knull{A}$ formed by the classes $[P_1],\ldots,[P_n]$, and in particular is again an integral matrix.

Similar situations may be encountered in
Lie theory, graph theory, knot theory, and singularity theory.
In the present paper we concentrate on the links to graph theory and to singularity theory. We show --- in the setting investigated in Sections~\ref{sect:singularities} and \ref{sect:stable} --- that the representation theory of a path algebra of a Dynkin quiver, or of a tubular canonical algebra, or of an extended canonical algebra yields a categorification of the \define{Milnor lattice}, a central tool for the analysis of isolated singularities. Compare~\cite{Ebeling:2007} for the relevant concepts.

The paper is organized as follows. In Section~\ref{sect:spectral_general} we review general properties of the spectrum and discuss the interplay between an algebra and its spectrum. In particular, we deal with the relationship to graph theory, in connection with path algebras of bipartite quivers. In Section~\ref{sect:important_classes} we review basic results concerning hereditary, canonical, extended canonical and supercanonical algebras. These classes of algebras are particularly suitable to spectral analysis. In Section~\ref{sect:further_prop} we consider the problem to characterize the derived class of an algebra by its spectral data. We discuss known and new examples of isospectral algebras. We consider the recently introduced concept of representability of self-reciprocal polynomials, relating certain Coxeter polynomials to spectra of graphs, thus extending the range of phenomena in graph theory related to representation theory.

In Section~\ref{sect:one_point} we study an important device to construct large algebras from smaller ones, the formation of one-point extensions, a process inverse to the formation of perpendicular categories with respect to exceptional modules. Explicit formulas are given expressing the corresponding changes of the Coxeter polynomials.  We further discuss properties of such extensions for the classes of hereditary and canonical algebras.

In Section~\ref{sect:accessible} we introduce a new class of algebras, called accessible algebras, which can be obtained by successive one-point extensions with exceptional modules, starting with the field $k$. Poset and tree algebras belong to this class, moreover the derived closure of the class of accessible algebras contains many of the algebras considered before. Despite the fact that also accessible algebras are not characterized by their Coxeter polynomials, we present a method to identify the derived class of an accessible algebra by spectral techniques.

In Section~\ref{sect:singularities} we present the (graded) singularities attached to a weighted projective line $\XX$, and discuss the shape of the triangulated category of the singularities  of $R$~(\cite{Buchweitz:1987}, \cite{Orlov:2005} and \cite{Krause:2005}), which largely depends on the sign of the Euler characteristic of $\XX$.
In the study of the derived category of finite dimensional modules over an extended canonical algebra, recently completed by the authors, an important tool was Orlov's theorem~\cite{Orlov:2005} to link this category to the triangulated category of singularities of a graded singularity, naturally associated to the context. This connection is the subject of sections~\ref{ssect:orlov} and~\ref{ssect:appl_orlov}.

In Section~\ref{sect:stable} we enlarge the range of hereditary representation theory by studying the stable category of vector bundles $\stabvect\XX$ on a weighted projective line.
This category inherits the structure of a triangulated category from the Frobenius category $\vect\XX$ of vector bundles on $\XX$, whose class of projective (=injective) objects is a naturally chosen class of line bundles on $\XX$. The stable category of vector bundles on $\XX$ is more accessible than the equivalent triangulated category of singularities studied in Section~\ref{sect:singularities}. For instance it allows an easy access to its Auslander-Reiten quiver. Based on spectral analysis, for the minimal wild weight type $(2,3,7)$, we offer a direct analysis of $\stabvect\XX$, thus                                                                                                                                                                                                                                                                                                                                                                                                                                                                                                                        bypassing Orlov's theorem. In the last Section~\ref{sect:comments} we offer a global view on the known classes of derived accessible algebras and discuss a number of open problems.

As standard references we mention \cite{Assem:Simson:Skowronski:2006}, \cite{Simson:Skowronski:2007a}, \cite{Simson:Skowronski:2007b}, \cite{Auslander:Reiten:Smalo:1995}, \cite{Ringel:1984} for algebras,  and \cite{Ebeling:2007} or singularities.

\section{The spectrum, general properties} \label{sect:spectral_general}
For the moment $k$ may be an arbitrary field. (Only later we will assume that $k$ is algebraically closed.) Let $A$ be a finite dimensional $k$-algebra. Let $S_1,\ldots,S_n$ be a complete system of pairwise nonisomorphic simple $A$-modules. Let $P_1,\ldots,P_n$ (respectively $I_1,\ldots,I_n$) denote complete systems of indecomposable projective (respectively injective) $A$-modules such that $P_i$ (resp. $I_i$) is the projective (resp. injective) hull of $S_i$. By Jordan-H\"{o}lder's theorem the classes $[S_1],\ldots,[S_n]$ of simple $A$-modules naturally form a basis of the Grothendieck group $\Knull{A}$ of the category $\mmod{A}$ of finite dimensional (right) $A$-modules with respect to short exact sequences.

If $A$ has finite global dimension then we dispose of a number of additional features:

(i) The Grothendieck group $\Knull{\Der{A}}$ of the bounded derived category $\Der{A}=\Der{\mmod{A}}$ of finitely generated $A$-modules is naturally isomorphic to $\Knull{A}$ by mapping a complex $C$ from $\Der{A}$ to the alternating sum of classes $[\homol{n}{C}]$.

(ii) The classes of indecomposable projectives (respectively injectives) form a basis of $\Knull{A}$.

(iii) The Auslander-Reiten transformation $\tau:\Der{A}\ra\Der{A}$ is a triangle equivalence, hence induces an automorphism $\cox{A}$ of $\Knull{A}$, the Coxeter transformation of $A$.

$(iv)$ $\Knull{A}$ is equipped with a (usually non-symmetric) bilinear homological form, given on classes of modules by
 $
 \eulerform{[X]}{[Y]}=\sum_{n\in\ZZ}(-1)^n\ddim{k}{\Ext{n}{}{X}{Y}}.
 $
This form is called the \define{Euler form}, it is non-degenerate.

$(v)$ Due to Serre duality of $\Der{A}$, Euler form and Coxeter transformation on $\Knull{A}$ are related by the formula $\eulerform{y}{x}=-\eulerform{x}{\cox{A}\,y} \textrm{ for all } x,y\in\Knull{A}.$
Hence $\Knull{A}$, equipped with the Euler form, is a bilinear lattice in the sense of \cite{Lenzing:1996}.

Fixing a basis $e_1,\ldots,e_n$ we may identify $\Knull{A}$ with $\ZZ^n$, with members written as column vectors. With the help of $C=(\eulerform{e_i}{e_j})$, called the \define{Cartan matrix} with respect to $e_1,\ldots,e_n$, we express the Euler form as
$\eulerform{x}{y}=\transpose{x}\,C\,y$.
The next assertion is an easy consequence of property $(iv)$ above.
\begin{lemma}
With the above conventions, the Coxeter transformation $\cox{A}$ is given by left matrix multiplication with the \define{Coxeter matrix} $\cox{A}=-\inverse{C}\transpose{C}$ with respect to $e_1,\ldots,e_n$.~\hfill\qed
\end{lemma}
While $C$ has a nonzero determinant, it may not be invertible over the integers, however the Coxeter matrix $\cox{A}$ is always an integral matrix (and of determinant one). This poses some restriction on the integral matrices qualifying as Cartan matrices. For the applications to follow, we mostly deal with the basis of  indecomposable projectives or simples, respectively. It will always be clear from the context, which of the two cases is considered.

The lemma implies that the \define{Coxeter polynomial} $\coxpol{A}$, that is, the characteristic polynomial of $\cox{A}$, is always monic, integral, and self-reciprocal. We recall that a polynomial $f=a_0+a_1x+\cdots+a_nx^n$ of degree $n$ is called \emph{self-reciprocal} (or symmetric, or palindromic) if $a_i=a_{n-i}$ for all $i=0,\ldots,n$. It is equivalent to request that $f(x)=x^nf(1/x)$. In this survey we will only consider self-reciprocal polynomials which are monic and integral.

By design, the Coxeter polynomial $\coxpol{A}$ reflects important homological properties of the algebra $A$. The exact nature of this relationship remains however mysterious, and only a few results of a general nature are known. The most notable instance is a result of Happel~\cite{Happel:trace}, expressing the degree one coefficient of $\coxpol{A}$ in terms of \define{Hochschild cohomology} $\hoch{i}{A}$. Note that $\hoch{i}{A}=0$ if $i$ is greater than the global dimension $\gldim{A}$ of $A$.
\begin{proposition}[Happel]\label{prop:Happel_trace}
The negative trace of the Coxeter transformation, hence the degree one coefficient of $\coxpol{A}$, equals  $\sum_{i=0}^\infty(-1)^i \ddim{k}{\hoch{i}{A}}$.~\hfill\qed
\end{proposition}

\subsection{Spectral radius one, periodicity} \label{ssect:periodicity}
If the spectrum of $A$ lies in the unit disk, then all roots of $\coxpol{A}$ lie on the unit circle, hence $A$ has spectral radius $\specrad{A}=1$. Clearly, for fixed degree there are only finitely many monic integral polynomials with this property.  Due to Kronecker's theorem, see \cite[Prop. 1.2.1]{Goodman:Harpe:Jones:1989}, these are easy to classify. We recall that the $n$-th \define{cyclotomic polynomial} $\Phi_n$ is the minimal polynomial of a primitive $n$-th root of unity over the rational number field $\QQ$. The polynomial $\Phi_n$ is monic integral of degree $\varphi(n)$, where $\varphi$ is Euler's totient function, see~\cite{Lang:2002}. The $\Phi_n$ can be recursively obtained from the formula $(x^n-1)=\prod_{d|n}\Phi_d$.

\begin{proposition}[Kronecker]
Let $f$ be a monic integral polynomial whose spectrum is contained in the unit disk. Then all roots of $f$ are roots of unity. Equivalently, $f$ factors into cyclotomic polynomials.~\hfill\qed
\end{proposition}
The following table displays the number of such polynomials $f$ for small degrees; $a(n)$ is the number of monic polynomials of degree $n$ of spectral radius $1$, $b(n)$ is the number of those which are additionally self-reciprocal and $c(n)$ is the number of those which are self-reciprocal and where $f(-1)$ is the square of an integer. (The reason to consider such polynomials will become clear later, see Section~\ref{ssect:triangular}).
\small
\begin{center}
\tabcolsep4pt
\begin{tabular}{|l|cccccccccccc|c|c|c|}\hline
$n$   &1& 2& 3& 4& 5& 6& 7& 8& 9& 10& 11& 12& 15& 20& 25 \\ \hline
$a(n)$&2& 6& 10& 24& 38& 78& 118& 224& 330& 584& 838& 1420& 4514& 30532& 152170\\
$b(n)$&1& 5& 5& 19& 19& 59& 59& 165& 165& 419& 419& 1001& 2257& 20399& 76085\\
$c(n)$&1& 3& 5& 12& 19& 34& 59& 99& 165& 244& 419& 598& 2257& 12526& 76085\\
 \hline
\end{tabular}
\end{center}
\normalsize
Indeed, there is an efficient algorithm to determine such polynomials of given degree $n$. The algorithm is based on a quadratic bound  $n\leq 4\,\totient(n)^2$ for $n$ in terms of Euler's totient $\totient(n)$, see \cite[p.\ 248]{Sierpinski:1988}. Note that there is no linear bound for $n$ in terms of $\totient(n)$.

Cyclotomic polynomials $\Phi_n$ and their products are a natural source for self-reciprocal polynomials. Clearly, $\Phi_1=x-1$ is not self-reciprocal, but all the remaining $\Phi_n$ (with $n\geq2$) are. Hence, exactly the polynomials $(x-1)^{2k}\prod_{n\geq2}\Phi_n^{e_n}$ with natural numbers $k,e_n$ are self-reciprocal of spectral radius one.

It is not a coincidence that in the above table we have $a(n)=b(n+1)$ for $n$ even and $a(n)=b(n)$ for $n$ odd. Indeed, if $f$ is self-reciprocal of odd degree then $f(-1)=0$ and hence $f=(x+1)g$, where $g$ is also self-reciprocal.

The following finite dimensional algebras are known to produce Coxeter polynomials of spectral radius one:
\begin{enumerate}
\item  hereditary algebras of finite or tame representation type, see Section~\ref{ssect:hereditary};
\item all canonical algebras, see Section~\ref{ssect:canonical};
\item (some) extended canonical algebras, see Section~\ref{ssect:ext_canonical};
\item generalizing (2), (some) algebras which are derived equivalent to categories of coherent sheaves.
\end{enumerate}
We put $v_n=1+x+x^2+\ldots+x^{n-1}$. Note that $v_n$ has degree $n-1$. There are several reasons for this choice: first of all $v_n(1)=n$, second this normalization yields convincing formulas for the Coxeter polynomials of canonical algebras and hereditary stars, third representing a Coxeter polynomial --- for spectral radius one --- as a rational function in the $v_n$'s relates to a Poincar\'{e} series, naturally attached to the setting, compare Section~\ref{sect:one_point}.
\small
\begin{center}
\renewcommand\arraystretch{1.4}
\begin{tabular}{|l|c|c|c|c|} \hline
  Dynkin  & star &$v$-factorization & cyclotomic & Coxeter \\
  type & symbol &  & factorization & number \\ \hline
  $\A{n}$& $[n]$ & $v_{n+1}$ & $\displaystyle\prod_{d|n,d>1}\Phi_d$ & $n+1$\\ \hline
  $\D{n}$&$[2,2,n-2]$& $\displaystyle\frac{v_2\,(v_2v_{n-2})}{(v_2v_{n-2})v_{n-1}}v_{2(n-1)}$ & $\displaystyle\Phi_2\,\prod_{\ofrac{d|2(n-1)}{d\neq 1,d\neq n-1}}\Phi_{d}$ & $2(n-1)$ \\\hline
  $\E{6}$ &$[2,3,3]$& $\displaystyle\frac{v_2 v_3 (v_3)}{(v_3)v_4v_6}v_{12}$ & $\Phi_3\Phi_{12}$& 12 \\ \hline
  $\E{7}$&$[2,3,4]$& $\displaystyle\frac{v_2v_3(v_4)}{(v_4)v_6 v_9}v_{18}$& $\Phi_2\Phi_{18}$& 18 \\\hline
  $\E{8}$& $[2,3,5]$ &$\displaystyle\frac{v_2v_3v_5}{v_6v_{10}v_{15}}v_{30}$&$\Phi_{30}$&30\\\hline
\end{tabular}
\end{center}
\normalsize
In the column `$v$-factorization', we have added some extra terms in the nominator and denominator which obviously cancel. The reason to complete the fraction in this way will become apparent in Section~\ref{ssect:graded_sing}.

Inspection of the table shows the following result:

\begin{proposition} \label{prop:dynkin}
Let $k$ be an algebraically closed field and $A$ be a connected, hereditary $k$-algebra which is representation-finite. Then the Coxeter polynomial $\coxpol{A}$ determines $A$ up to derived equivalence.~\hfill \qed
\end{proposition}
\subsection{Triangular algebras} \label{ssect:triangular}
Nearly all algebras considered in this survey are triangular. By definition, a finite dimensional algebra is called \define{triangular} if it has triangular matrix shape
\small
$$
\left[
\begin{array}{ccccc}
A_1 & M_{12}&\cdots &M_{1n}\\
0   &A_2    &\cdots &M_{2n}\\
    & & \ddots  & \vdots\\
0   & 0     &\cdots&A_{n}
\end{array}
\right]
$$\normalsize
where the diagonal entries $A_i$ are skew-fields and the off-diagonal entries $M_{ij}$, $j>i$, are $A_i,A_j$-bimodules. Each triangular algebra has finite global dimension.
\begin{proposition}
  Let $A$ be a triangular algebra over an algebraically closed field $k$. Then  $\coxpol{A}(-1)$ is the square of an integer.
\end{proposition}
\begin{proof}
  Let $C$ be the Cartan matrix of $A$ with respect to the basis of indecomposable projectives. Since $A$ is triangular and $k$ is algebraically closed, we get $\det C=1$, yielding
  $$
  \coxpol{A}=\det\left( x I + C^{-1}C^t \right)=\det\left(C^{-1} \right)\cdot\det\left(x C+C^t\right)=\det\left(C^t+xC \right).
  $$
  Hence $\coxpol{A}(-1)$ is the determinant of the skew-symmetric matrix $S=C^t-C$. Using the skew-normal form of $S$, see \cite[Theorem IV.1]{Newman:1972},
  we obtain  $S'=U^t S U$ for some $U\in\GL{n}\ZZ$, where
  $S'$ is a
  block-diagonal matrix whose first block is
  the zero matrix of a certain size and where the remaining blocks have the shape
  \footnotesize$
  \left[
    \begin{matrix}
      0 & m_i \\ -m_i & 0
    \end{matrix}
  \right]
  $\normalsize with integers $m_i$. The claim follows.
\end{proof}
\noindent\emph{Which self-reciprocal polynomials of spectral radius one are Coxeter polynomials?} The answer is not known. If arbitrary base fields are allowed, we conjecture that all self-reciprocal polynomials are realizable as Coxeter polynomials of triangular algebras. Restricting to algebraically closed fields, already the request that $f(-1)$ is a square discards many self-reciprocal polynomials, for instance the cyclotomic polynomials $\Phi_4$, $\Phi_6$, $\Phi_8$, $\Phi_{10}$. Moreover, the polynomial $f=x^3+1$, which is the Coxeter polynomial of the non simply-laced Dynkin diagram $\B{3}$,  does not appear as the Coxeter polynomial of a triangular algebra over an algebraically closed field, despite of the fact that $f(-1)=0$ is a square. Indeed, the Cartan matrix
\small
$$
\left[
\begin{array}{ccc}
1&a&b\\
0&1&c\\
0&0&1\\
\end{array}
\right]
$$\normalsize
yields the Coxeter polynomial $f=x^3+\alpha x^2 +\alpha x+1$, where $\alpha=abc-a^2-b^2-c^2+3$. The equation $a^2+b^2+c^2-abc=3$ of Hurwitz-Markov type does not have an integral solution. (Use that reduction modulo $3$ only yields the trivial solution in $\FF_3$.)

\subsection{Relationship with graph theory} \label{ssect:graph}

Given a (non-oriented) graph $\De$, its \define{characteristic polynomial} $\charpol{\De}$ is defined as the characteristic polynomial of the adjacency matrix $M_\De$ of $\De$.
Observe that, since $M_\De$ is symmetric, all its eigenvalues are real numbers. For general results on graph theory and spectra of graphs see \cite{Cvetkovic:Doob:Sachs:1979} and \cite{Cvetkovic:Rowlinson:Simic:1997}.

There are important interactions between the theory of graph spectra and the representation theory of algebras, due to the fact that
if $C$ is the Cartan matrix of $A= k[\vDe]$, then $M_\De$ is determined by the symmetrization $C+C^t$ of $C$, since $M_\De =C+C^t-2I$. We shall see that information on the spectra of $M_\De$ provides fundamental insights into the spectral analysis of the Coxeter matrix $\Phi_A$ and the structure of the algebra $A$.

A fundamental fact for a hereditary algebra $A=k[\vDe]$, when $\vDe$ is a
\define{bipartite quiver}, that is, every vertex is a sink or source, is that
$\Spec{\cox{A}}\subset \SS^1\cup \RR^+$, see Section~\ref{ssect:representability}. This was shown by A'Campo
\cite{ACampo:1973} as a consequence of the following important identity.

\begin{proposition}[A'Campo] \label{prop:ACampo}
Let $A=k[\vDe]$ be a hereditary algebra with $\vDe$ a
\emph{bipartite quiver\/} without oriented cycles. Then
$\coxpol{A}(x^2)=x^n \charpol\De(x+x^{-1})$, where $n$ is the number of vertices of $\vDe$
and $\charpol\De$ is the characteristic polynomial of the
underlying graph $\De$ of $\vDe$
\end{proposition}
\noindent\emph{Proof.}
Since $\vDe$ is bipartite, we may assume that the first $m$ vertices are sources and the last $n-m$ vertices are sinks. Then the adjacency matrix $A$ of $\De$ and the Cartan matrix $C$ of $A$, in the basis of simple modules, take the form: $A= N + N^t$, $C= I_n - N$, where
$$N= \begin{pmatrix} 0& D \cr 0 & 0 \cr \end{pmatrix}$$
for certain $m \times m$-matrix $D$. Since $N^2 = 0$, then $C^{-1}=I_n + N$. Therefore
\begin{align*}
  \det(x^2I_n-\cox{A}) &=& \det(x^2I_n +(I_n-N)(I_n+N)^t) \det(I_n-N^t)\\
  =\det(x^2I_n-x^2N^t+(I_n-N)) &=& x^n \det((x+x^{-1})I_n-xN^t-x^{-1}N)\\
  &=&x^n \det((x+x^{-1})I_n-A).\end{align*}~\vskip-0.9cm\hfill\qed

We shall come back to this \define{representability} property of $\coxpol{A}$.

The above result is important since it makes the spectral analysis of bipartite quivers and their underlying graphs almost equivalent. Note, however, that the representation theoretic context is much richer, given the categorical context behind the spectral analysis of quivers. The representation theory of bipartite quivers may thus be seen as a categorification of the class of graphs, allowing a bipartite structure.

\subsubsection*{Constructions in graph theory} Several simple constructions in graph theory provide tools to obtain in practice the characteristic polynomial of a graph. We recall two of them (see \cite{Cvetkovic:Doob:Sachs:1979} for related results):

(a) Assume that $a$ is a vertex in the graph $\De$ with a unique neighbor $b$ and $\De'$(resp. $\De''$) is the full subgraph of $\De$ with vertices $\De_0\setminus \{a\}$ (resp. $\De_0\setminus \{a,b\}$), then
$$\charpol\De=x \charpol{\De'}- \charpol{\De''}$$

(b) Let $\De_i$ be the graph obtained by deleting the vertex $i$ in $\De$. Then the first derivative of $\charpol{\De}$ is given by
$$\charpol{\De}' = \sum_i \charpol{\De_i}$$
The above formulas can be used inductively to calculate the characteristic polynomial of trees and other graphs. They immediately imply the following result that will be used often to calculate Coxeter polynomials of algebras.

\begin{proposition}
Let $A=k[\vDe]$ be a bipartite hereditary algebra. The following holds:

(i) Let $a$ be a vertex in the graph $\De$ with a unique neighbor $b$. Consider the algebras $B$ and $C$ obtained as quotients of $A$ modulo the ideal generated by the vertices $a$ and $a,b$, respectively. Then
$$\coxpol{A}=(x+1) \coxpol{B}- x \coxpol{C}.$$

(ii) The first derivative of the Coxeter polynomial satisfies:
$$2x \coxpol{A}' =n \coxpol{A} + (x-1) \sum_{i}\coxpol{A^{(i)}}$$
where $A^{(i)}=k[\vDe \setminus{\{i\}}]$ is an algebra obtained from $A$ by `killing' a vertex $i$.
\end{proposition}
\begin{proof}
Use the corresponding results for graphs and A'Campo's formula for the algebras $A$ and its quotients $A^{(i)}$.
\end{proof}
Later part $(i)$ of the proposition will be extended to triangular algebras, see Proposition~\ref{prop:coxpol_onepoint}.
\section{Important classes of algebras} \label{sect:important_classes}
In this section we give the definitions and main properties of such classes of finite dimensional algebras where information on their spectral properties is available. It is no accident that these algebras will reappear in Section~\ref{sect:accessible}, where we are going to describe a powerful method to decide on the (derived) shape of an algebra through spectral analysis.

\subsection{Hereditary algebras} \label{ssect:hereditary}

Let $A$ be a finite dimensional $k$-algebra with $k$ an algebraically
closed field. For simplicity we assume $A=k[\vDe]/I$ for a quiver $\vDe$ without
oriented cycles and $I$ an ideal of the path algebra.
The following facts about the Coxeter transformation $\cox{A}$ of $A$ are fundamental:

(i) Let $S_1,\ldots,S_n$ be a complete system of pairwise
non-isomorphic simple $A$-modules, $P_1,\ldots,P_n$ the corresponding
projective covers and $I_1,\ldots,I_n$ the injective envelopes. Then
$\cox{A}$ is the automorphism of $\Knull{A}$ defined by $\cox{A}[P_i]=-[I_i]$, where $[X]$
denotes the class of a module $X$ in $K_0(A)$.

(ii) For a hereditary algebra $A=k[\vDe]$, the \define{spectral radius}
$\specrad{A}=\specrad{{\cox{A}}}$ determines the representation type of $A$ in the
following manner:

$\bullet$ $A$ is representation-finite if $1=\specrad{A}$
is not a root of the Coxeter polynomial $\coxpol{A}$.

$\bullet$ $A$ is tame if $1=\specrad{A}\in\Roots{\coxpol{A}}$.

$\bullet$ $A$ is wild if $1<\specrad{A}$. Moreover, if $A$ is wild connected,
Ringel \cite{Ringel:1994} shows that the spectral radius $\specrad{A}$ is a simple root of
$\coxpol{A}$. Then Perron-Frobenius theory yields a vector $y^+\in
K_0(A)\otimes_{\ZZ}\RR$ with positive coordinates such
that $\cox{A}y^+ =\specrad{A}y^+$.  Since $\coxpol{A}$ is
self-reciprocal, there is a vector $y^-\in
K_0(A)\otimes_{\ZZ}\RR$ with positive coordinates such
that $\cox{A}y^- =\specrad{A}^{-1}y^-$. The vectors $y^+,y^-$ play an
important role in the representation theory of $A=k[\vDe]$, see \cite{Dlab:Ringel:1981},
\cite{Pena:Takane:1990}. Namely, for an indecomposable $A$-module $X$:

{\rm (a)} $X$ is a preprojective $A$-module (that is,
$\tau^m_AX$ is projective for some $m\ge 0$) if and only if $\eulerform{y^-}{[X]}_A<0$

{\rm (b)} $X$ is a preinjective $A$-module (that is,
$\tau^{-m}_AX$ is injective for some $m\ge 0$) if and only if
$\eulerform{[X]}{y^+}_A<0$.

{\rm (c)} $X$ is regular (that is, $X$ is not preprojective or
preinjective) if and only if $\eulerform{y^-}{[X]} >0$ and $\eulerform{[X]}{y^+}>0$.

{\rm (d)} If $X$ is preprojective or regular, then
$\!\lim\limits_{n\to \infty}\frac{1}{\specrad{A}^n}[\tau^{-n}_AX]\!=\!\lambda^-_Xy^-$, for some
$\lambda^-_X\!>\!0$.

{\rm (e)} If $X$ is preinjective or regular, then
$\lim\limits_{n\to \infty}{\frac{1}{\specrad{A}^n}[\tau^n_AX]}=\lambda^+_Xy^+$, for some $\lambda^+_X>0$.

The above criteria for indecomposable modules over hereditary algebras provide a shortcut into deep results in the representation theory, an example:
\begin{theorem}
Let $X, Y$ be indecomposable regular modules over a wild hereditary algebra $A$. Then there is a number $N$ such that for every $m > N$ we have:

(a) \emph{\cite{Baer:1986}}: $\Hom{A}{X}{\tau^m Y} \ne 0$;

(b) \emph{\cite{Kerner:1995}}: $\Hom{A}{X}{\tau^{-m} Y} = 0$.

In particular, given two regular components $C_1,C_2$ we have $\Hom{A}{C_1}{C_2} \ne 0$. Inside a regular component, most of the morphisms go in the direction opposite to the arrows.
\end{theorem}

\begin{proof}
(a): For some number $\lambda_Y^+ > 0$, we have
$\!\lim\limits_{n\to \infty}{\frac{1}{\specrad{A}^n}[\tau^m_AY]}\!=\!\lambda^+_Yy^+$
Therefore $ 0 < \eulerform{[X]}{y^+}_A =\!\lim\limits_{n\to \infty}\frac{1}{\specrad{A}^n}
\eulerform{[X]}{[\tau^m_AY]}_A$.
\end{proof}

\subsubsection*{Explicit formulas, special values}
We are discussing various instances where an explicit formula for the Coxeter polynomial is known.
\subsubsection*{star quivers}
Let $A$ be the path algebra of a hereditary star $[p_1,\ldots,p_t]$
with respect to the standard orientation, see
 $$
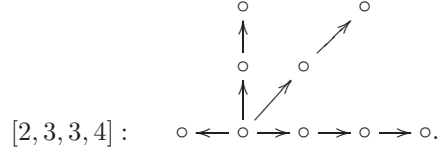
\def\c{\circ}
\xymatrix@C12pt@R12pt{
        &&\c                          &        &\c  &\\
        &&\c\ar[u]                    &\c\ar[ru]&     &\\
[2,3,3,4]:&\c      &\c\ar[l]\ar[u]\ar[r]\ar[ru]&\c\ar[r]&\c\ar[r]&\c.\\
}
$$
Since the Coxeter polynomial $\coxpol{A}$ does not depend on the orientation of $A$ we will denote it by $\coxpol{[p_1,\ldots,p_t]}$. It follows from \cite[prop.~9.1]{Lenzing:Pena:1997} or \cite{Boldt:1995} that
\begin{equation}\label{eq:star_formula}
\coxpol{[p_1,\ldots,p_t]}=\prod_{i=1}^{t}v_{p_i}\left((x+1)-x\,\sum_{i^=1}^t \frac{v_{p_i-1}}{v_{p_i}} \right).
\end{equation}
In particular, we have an explicit formula for the sum of coefficients of $\chi=\coxpol{[p_1,\ldots,p_t]}$ as follows:
$$
\chi(1)=\prod_{i=1}^t p_i\left( 2-\sum_{i=1}^t (1-\frac{1}{p_i})\right).
$$
This special value of $\chi$ has a specific mathematical meaning:
up to the factor $\prod_{i=1}^t p_i$ this is just the orbifold-Euler characteristic of a weighted projective line $\XX$ of weight type $(p_1,\ldots,p_t)$. Moreover,
\begin{enumerate}
\item $\chi(1)>0$ if and only if the star $[p_1,\ldots,p_t]$ is of Dynkin type, correspondingly the algebra $A$ is representation-finite.
\item $\chi(1)=0$ if and only if the star $[p_1,\ldots,p_t]$ is of extended Dynkin type, correspondingly the algebra $A$ is of tame (domestic) type.
\item $\chi(1)<0$ if and only if $[p_1,\ldots,p_t]$ is not Dynkin or extended Dynkin, correspondingly the algebra $A$ is of wild representation type.
\end{enumerate}
The above deals with all the Dynkin types and with the extended Dynkin diagrams of type $\Dext{n}$, $n\geq 4$, and $\Eext{n}$, $n=6,7,8$. To complete the picture, we also consider the extended Dynkin quivers of type $\Aext{n}$ ($n\geq 2$) restricting, of course, to quivers without oriented cycles. Here, the Coxeter polynomial depends on the orientation: If $p$ (resp.\ $q$) denotes the number of arrows in clockwise (resp.\ anticlockwise) orientation ($p,q\geq1, p+q=n+1$), that is, the quiver has type $\Aext{p,q}$, the Coxeter polynomial $\chi$ is given by
$$
\coxpol{(p,q)}=(x-1)^2\,v_{p}v_{q}.
$$
Hence $\chi(1)=0$, fitting into the above picture.

The next table displays the $v$-factorization of extended Dynkin quivers
\bigskip
\begin{center}
\renewcommand\arraystretch{1.3}
\begin{tabular}{|c|c|c|c|}\hline
extended Dynkin type & star symbol & weight symbol & Coxeter polynomial \\ \hline
$\Aext{p,q}$           &     ---     & $(p,q)$         & $(x-1)^2v_p\,v_q$ \\ \hline
$\Dext{n}$, $n\geq4$   & [2,2,n-2]   & $(2,2,n-2)$   & $(x-1)^2v_2^2v_{n-2}$\\\hline
$\Eext{6}$& $[3,3,3]$ & $(2,3,3)$    & $(x-1)^2 v_2 v_3^2$ \\\hline
$\Eext{7}$& $[2,4,4]$ & $(2,3,4)$    & $(x-1)^2 v_2v_3v_4$ \\\hline
$\Eext{8}$& $[2,3,6]$ & $(2,3,5)$    & $(x-1)^2 v_2v_3v_5$\\\hline
\end{tabular}
\end{center}

\begin{remark}
As is shown by the above table, proposition~\ref{prop:dynkin} extends to the tame hereditary case. That is, the Coxeter polynomial of a connected, tame hereditary $k$-algebra $A$ ($k$ algebraically closed) determines the algebra $A$ up to derived equivalence. This is no longer true for wild hereditary algebras, not even for trees, as will be shown in Section~\ref{ssect:isospectral}.
\end{remark}

\subsection{Canonical algebras} \label{ssect:canonical}
A canonical algebra $\La=\Lapla$ is determined by a weight sequence $\pP=(p_1,\ldots,p_t)$ of $t\geq2$ integers $p_i\geq2$ and a parameter sequence $\lala=(\la_3,\ldots,\la_t)$ consisting of $t-2$ pairwise distinct non-zero scalars  from the base field $k$. (We may assume $\la_3=1$ such that for $t\le 3 $ no parameters occur).  Then the algebra $\La=\Lapla$ is defined by the quiver
$$
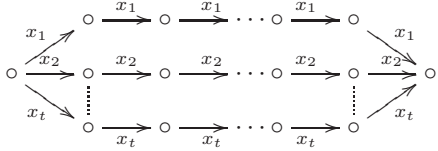
\def\c{\circ}
\xymatrix@C18pt@R10pt{
                                        &\c\ar[r]^{x_1}          & \c\ar[r]^{x_1}         &\cdots\c\ar[r]^{x_1}&\c\ar[dr]^{x_1}&\\
\c\ar[r]^{x_2}\ar[ru]^{x_1}\ar[rd]_{x_t}&\c\ar[r]^{x_2}\ar@{.}[d]&\c\ar[r]^{x_2}          &\cdots\c\ar[r]^{x_2}&\c\ar[r]^{x_2}\ar@{.}[d] &\c    \\
                                        &\c\ar[r]_{x_t}          &\c\ar[r]_{x_t}          &\cdots\c\ar[r]_{x_t}&\c\ar[ru]_{x_t}&    \\
}
$$
satisfying the $t-2$ equations:
$$x_i^{p_i}=x_1^{p_1}-\lambda_i x_2^{p_2},\qquad i=3,\ldots,t.$$
For more than two weights, canonical algebras are not hereditary. Instead, their representation theory is determined by a hereditary category, the category $\coh\XX$ of coherent sheaves on a \define{weighted projective line} $\XX=\Xpla$, naturally attached to $\La$, see~\cite{Geigle:Lenzing:1987}.
\begin{proposition}
  Let $\La=\Lapla$ be a canonical algebra. Then $\La$ is the endomorphism ring of a tilting object in the category $\coh\XX$ of coherent sheaves on the weighted projective line $\XX=\Xpla$. The category $\coh\XX$ is hereditary and satisfies Serre duality in the form $\dual{\Ext1{}{X}{Y}}=\Hom{}{Y}{\tau X}$ for a self-equivalence $\tau$ which serves as the Auslander-Reiten translation.~\hfill\qed
\end{proposition}
Canonical algebras were introduced by Ringel~\cite{Ringel:1984}. They form a key class to study important features of representation theory. In the form of tubular canonical algebras they provide the standard examples of tame algebras of linear growth. Up to tilting canonical algebras are characterized as the connected $k$-algebras with a separating exact subcategory or a separating tubular one-parameter family (see \cite{Lenzing:Pena:1999} and \cite{Skowronski:1996}).
That is, the module category $\mmod\La$ accepts a \define{separating tubular family}
$\mathcal{T}=(T_\lambda)_{\lambda \in
P_1k}$, where $T_\lambda$ is a homogeneous tube for all $\lambda$
with the exception of $t$ tubes $T_{\lambda_1},\ldots,T_{\lambda_t}$
with $T_{\lambda_i}$ of rank $p_i$ ($1\le i\le t$).

Canonical algebras constitute an instance, where the explicit form of the Coxeter polynomial is known, see~\cite{Lenzing:Pena:1997} or \cite{Lenzing:1996}.
\begin{proposition}
Let $\La$ be a canonical algebra with weight and parameter data $(\pP,\lala)$. Then the Coxeter polynomial of $\La$ is given by
$$
\coxpol{\La}=(x-1)^2\prod_{i=1}^t v_{p_i}.
$$~
\vskip-1cm\hfill\qed
\end{proposition}
The Coxeter polynomial therefore only depends on the weight sequence $\pP$. Conversely, the Coxeter polynomial determines the weight sequence --- up to ordering.

A finite dimensional algebra isomorphic to the endomorphism algebra of a tilting object in a (connected) hereditary abelian Hom-finite $k$-category $\Hh$ is called a \define{quasi-tilted algebra}. By a result of Happel~\cite{Happel:2001} each such category is derived equivalent to the module category $\mmod{H}$ over a hereditary algebra or to the category $\coh\XX$ of coherent sheaves on a weighted projective line. The Coxeter polynomials of quasi-tilted algebras are therefore the Coxeter polynomials of hereditary or canonical algebras.

In Section~\ref{sect:singularities} we investigate a class of graded singularities naturally attached to weighted projective lines or canonical algebras. There we will provide more details on the hereditary category $\coh\XX$.

\subsection{Supercanonical algebras}\label{ssect:supercanonical}
Following \cite{Lenzing:Pena:2004}, see also \cite{Skowronski:2001},
a supercanonical algebra is defined as follows:
The \emph{double cone} $\hat{S}$ of a finite poset $S$ is the poset obtained from $S$ by adjoining a unique source $\al$ and a unique sink $\om$, like in the picture
$$
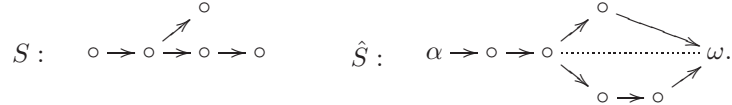
\def\c{\circ}
\def\b{\bullet}
\xymatrix@C10pt@R4pt{
   &         &                &\c      &  & &        &        &        &                &\c\ar[rrd]&         &  \\
S: & \c\ar[r]&\c\ar[ru]\ar[r] &\c\ar[r]&\c& &\hat{S}: &\al\ar[r]&\c\ar[r]&\c\ar[ur]\ar[dr]\ar@{.}[rrr]&          &         &\om.\\
   &         &                &        &  & &        &        &        &                &\c\ar[r]  &\c\ar[ur]&  \\
}
$$
Due to commutativity there is a unique path $\ka_S$ from $\al$ to $\om$ in $\hat{S}$.
Let now $S_1,\ldots,S_t$ be  finite posets, $t\ge 2$ and $\lambda_3,\ldots,\lambda_t$
pairwise different elements from  $ k\setminus \{0\}$. The \emph{supercanonical algebra}
$A=A(S_1,\ldots,S_t;\lambda_3,\ldots,\lambda_t)$ is obtained from the fully commutative quivers $\hat{S}_1,\ldots,\hat{S}_t$ by identification of the sources and sinks, respectively, and requesting additionally the
$t-2$ relations $\kappa_i=\kappa_2-\lambda_i\kappa_1$, $i=3,\ldots,t$, where $\ka_i=\ka_{S_i}$. The next figure yields an example of a supercanonical algebra with three arms
$$
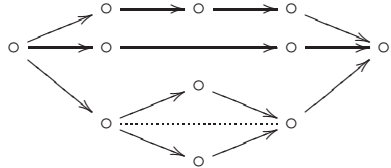
\def\c{\circ}
\xymatrix@R4pt{
                           &\c\ar[r]          &\c\ar[r]          &\c\ar[rd]        &  \\
\c\ar[r]\ar[ur]\ar[ddr]    &\c\ar[rr]          &                  &\c\ar[r]         &\c\\
                           &                  &\c\ar[dr]         &                 &  \\
                           &\c\ar[ur]\ar[dr]\ar@{.}[rr]&                  &\c\ar[uur]       &  \\
                           &                  &\c\ar[ur]         &                 &  \\
}
$$
where we assume that $\ka_3=\ka_2-\ka_1$.
In case, $S_1,\ldots,S_t$ are linear quivers
$S_i=[p_i]\colon 1\to 2\to \cdots \to p_{i-1}$,
the algebra $A(S_1,\ldots,S_t;\lambda_1,\ldots,\lambda_t)$
 is just the canonical algebra $\La(p_1,\ldots,p_t;\lambda_3,\ldots,\lambda_t)$.

Returning to the general case, a simple calculation shows that
$$\charpol{A}=(x-1)^2\prod^t_{i=1}\charpol{S_i}$$
where $\coxpol{S_i}$ is the Coxeter polynomial of the poset algebra
$S_i$, $i=1,\ldots,t$.

Supercanonical algebras form a natural generalization of canonical algebras,  since they arise from a canonical algebra $\La$ by replacing the linear arms of $\La$
by finite posets, that is, fully commutative quivers. In this paper we are mainly interested in \define{supercanonical algebras of restricted type}, where in addition to linear arms only posets of the form
 $$
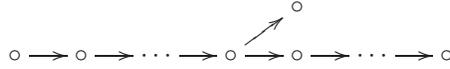
\def\c{\circ}
\xymatrix@C14pt@R8pt{
       &       &             &                   & \c         &         &    \\
\c\ar[r]&\c\ar[r]&\cdots\ar[r] &\c\ar[r]\ar[ru] & \c\ar[r] &\cdots\ar[r]&\c
}
$$
will appear. There are many reasons to make supercanonical algebras an interesting class:

$(1)$ The K-theory of supercanonical algebras is very similar to the K-theory of canonical algebras. There is a Riemann-Roch theorem, moreover one has an explicit formula for the Coxeter polynomial in terms of the characteristic polynomials of the posets $S_i$.

$(2)$ One-point extensions of canonical algebras with exceptional modules of derived finite length are supercanonical, see Theorem~\ref{thm:1pt-canonical}.

\subsection{Extended canonical algebras} \label{ssect:ext_canonical}
For an algebra $A$ and a right $A$-module $M$ we call $A[M]=$\small$\left[\begin{array}{cc}
A & 0\\ M &k
\end{array}\right]$\normalsize the \define{one-point extension} of $A$ by $M$. In Section~\ref{sect:one_point}
we will study such algebras in greater detail. Let $\La=\Lapla$ be a canonical
algebra. In \cite{Lenzing:Pena:2006} the authors proposed to study an interesting
class of algebras according with the following result:

\begin{proposition} \label{prop:ext_canonical}
The derived equivalence class of the one-point extension of a canonical algebra $\La$ by
an indecomposable projective or injective module is independent of
the particular choice of the module.
\end{proposition}

\begin{proof}
In the derived category $\Der\La=\Der{\coh\XX}$ an indecomposable projective or injective module $P$ is, up to translation, a line bundle from $\coh\XX$. For any two such objects $P_1$ and $P_2$ there hence exists a self-equivalence of $\Der{\coh\XX}$ sending
$P_1$ to $P_2$, see \cite[prop.~2.1]{Geigle:Lenzing:1987}. The assertion now
follows from Proposition~\ref{prop:one-point}.
\end{proof}

We call an algebra of the form $\La[P]$, with $P$ indecomposable
projective, an \define{extended canonical algebra\/} of type
$\spitz{\pP,\lala}$ or just of type $\spitz{p_1,\ldots,p_t}$. We use the notation $\Lahat\spitz{\pP,\lala}$ if $P$ is the indecomposable projective of the sink vertex of the quiver of $\La$. Note that  $\Lahat=\Lahat\spitz{\pP,\lala}$ is given by the quiver
$$
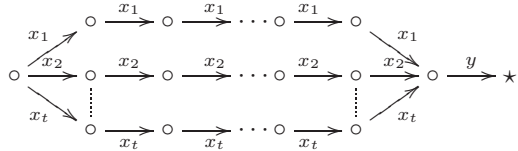
\def\c{\circ}
\xymatrix@C18pt@R10pt{
                                        &\c\ar[r]^{x_1}          & \c\ar[r]^{x_1}         &\cdots\c\ar[r]^{x_1}&\c\ar[dr]^{x_1}&         &\\
\c\ar[r]^{x_2}\ar[ru]^{x_1}\ar[rd]_{x_t}&\c\ar[r]^{x_2}\ar@{.}[d]&\c\ar[r]^{x_2}          &\cdots\c\ar[r]^{x_2}&\c\ar[r]^{x_2}\ar@{.}[d] &\c\ar[r]^y &\star\\
                                        &\c\ar[r]_{x_t}          &\c\ar[r]_{x_t}          &\cdots\c\ar[r]_{x_t}&\c\ar[ru]_{x_t}    &         & \\
}
$$
keeping all the relations for the canonical algebra $\La$ and not introducing new ones. In particular, there are no relations involving the new vertex $\star$.
For the Coxeter polynomial $\coxpol{\Lahat}$, which only depends on the
numbers $p_1,\ldots,p_t$, we write $\extcoxpol{\spitz{p_1,\ldots,p_t}}$ or $\coxpol{\spitz{p_1,\ldots,p_t}}$.
By Proposition~\ref{prop:onepoint_coxpol} we immediately get:
\begin{coro}
The extended canonical algebra $\Lahat$ of type $\spitz{\pP,\lala}$ has Coxeter
polynomial
$$\extcoxpol{\spitz{p_1,\ldots,p_t}}=(x+1)(x-1)^2\prod^t_{i=1}v_{p_i}
-x\coxpol{[p_1,\ldots,p_t]}.$$~\vskip-1.2cm\hfill\qed
\end{coro}
Recall that the explicit form of the Coxeter polynomial for the hereditary star $[p_1,\ldots,p_t]$ is given in Section~\ref{ssect:hereditary}.

The algebra $\hat{D}=\hat{D}\spitz{\pla}$  given in terms of the
quiver
$$\def\bu{\ci}
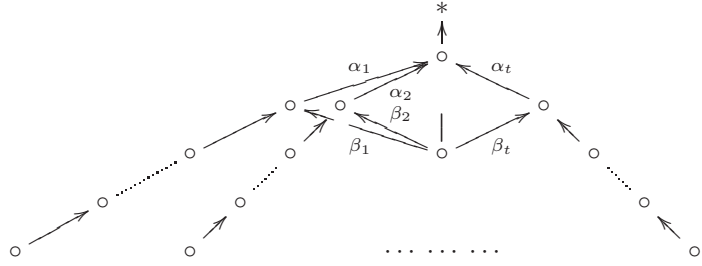
\xymatrix@R8pt@C8pt{
 &&&           &              &                &                   &                  &     \ast                            &                  &          &          &\\
 &&&           &              &                &                   &                  &     \bu\ar[u]                                         &                  &          &          &\\
 &&&           &              &                &\bu\ar[urr]^{\al_1}&\bu\ar[ur]_{\al_2}&                                                 &\bu\ar[lu]_{\al_t}&          &          &\\
 &&&           & \bu\ar[urr]   &                &\bu\ar[ur]         &                  &\bu\ar@{--}[u]\ar@{--}[u]\ar[ur]_{\be_t}\ar[llu]^{\be_1}|!{[u];[lllll]}\hole\ar[lu]_{\be_2}&                  &\bu\ar[lu]&          &\\
 &&\bu\ar@{.}[urr]&           &              &\bu\ar@{.}[ur]  &                   &                  &                                        &                  &          &\bu\ar@{.}[lu]&\\
 \bu\ar[urr]&&&   &\bu\ar[ur]    &                &                   &                  &                                         \cdots\;\cdots\;\cdots        &                  &         &          &\bu\ar[lu]\\
}$$ with the two relations $\sum_{i=2}^t
\al_i\be_i=0$ and $\al_1\be_1=\sum_{i=3}^t
\la_i\,\al_i\be_i$ is called the \define{Coxeter-Dynkin algebra of extended canonical type}. For the next statement we refer to~\cite{Lenzing:Pena:2006}.

\begin{proposition}
   A Coxeter-Dynkin algebra $\hat{D}$ and an extended canonical type $\Lahat$ of the same type $\spitz{\pP,\lala}$ are derived equivalent.~\hfill\qed
\end{proposition}

\subsubsection*{The derived category of an extended canonical algebra}
The structure of the bounded derived category of an extended canonical algebra
$\Lahat=\Lahat\spitz{\pP,\lala}$ sensibly depends on the sign of the Euler
characteristic $\eulerchar{\XX}=2-\sum_{i=1}^t(1-1/p_i)$ of the weighted projective line $\XX$
associated to $\La$.
Following \cite{Lenzing:Pena:2006},  the description of the derived category of an extended canonical algebra yields an interesting trichotomy. Before proceeding the reader is advised to read Section~\ref{ssect:one_point_perp}.

Assume $\Tt$ is a triangulated category which is algebraic in the sense of Keller~\cite{Keller:2006}. All triangulated categories appearing in this survey are algebraic. An exceptional object $E$ is called
\define{special in $\Tt$} if the left (resp.\ right) perpendicular category $\lperp{E}$ (resp.\ $\rperp{E}$) is equivalent to
$\Der{\coh\XX}$ for some weighted projective line $\XX$ and,
moreover, the left adjoint $\ell$ (resp.\ right adjoint $r$) to inclusion maps $E$ to a line bundle in
$\coh\XX$. The next proposition, taken from \cite{Lenzing:Pena:2006} is the key tool to determine the shape of the derived category of an extended canonical algebra.
\begin{proposition}\label{prop:special}
Let $\Tt$ be a triangulated category having an exceptional object
$E$ that is special in $\Tt$. Then there exists a tilting object
$\bar{T}$ of $\Tt$ whose endomorphism ring is an extended canonical
algebra. Further for $A=\End{}{\bar{T}}$ the categories $\Tt$ and $\Der{\mmod{A}}$ are equivalent as triangulated categories.
\end{proposition}
\begin{proof}
By \cite{Geigle:Lenzing:1987} the line bundle $rE$ of $\coh\XX$ extends to a tilting bundle $T$ in $\coh\XX$. Now $\bar{T}=T\oplus E$ is a tilting object in $\Tt$, whose endomorphism ring is an extended canonical algebra $\Lahat$. The claim follows.
\end{proof}

\subsubsection*{Positive Euler characteristic: the domestic case} Consider a canonical algebra $\La=\Lapla$ of domestic type,
that is, $\eulerchar\XX>0$. We thus can assume $t=3$, allowing some weights to be $1$, if necessary. Let $\Delta
=[p_1,p_2,p_3]$ be the corresponding Dynkin star and $\tilde{\Delta}$ be the corresponding
extended Dynkin diagram. Then $\tilde{\Delta}$ admits a unique
positive \define{additive function} $\lambda$ assuming the value $1$
in some vertex (called \define{extension vertex}). Additivity of $\la$ means that for any vertex $u$ one has $2\la(u)=\sum_{v }a_{uv}\la(v)$, where $v$ runs through all vertices and $a_{uv}$ denotes the number of vertices between $u$ and $v$.
The \define{double extended graph} of type $\Delta$, denoted by
$\tilde{\tilde{\Delta}}$, is the graph arising from $\tilde{\Delta}$
by adjoining a new edge in an extension vertex. We illustrate this for the case $\De=[2,3,3]$, where the diagram in the middle gives the values of the additive function.
\scriptsize$$
\def\c{\circ}
\begin{array}{ccc}
\xymatrix@C6pt@R6pt{
 &           &             &                       &             &   \\
 &           &             & \c                    &             &   \\
\De:&\c\ar@{-}[r]& \c\ar@{-}[r]& \c\ar@{-}[u]\ar@{-}[r]& \c\ar@{-}[r]& \c\\
} &
\xymatrix@C6pt@R6pt{
 &           &             &  1                    &             &   \\
 &           &             &  2 \ar@{-}[u]            &             &   \\
\tilde\De:&1\ar@{-}[r]& 2\ar@{-}[r]& 3\ar@{-}[u]\ar@{-}[r]& 2\ar@{-}[r]& 1\\
}
&
\xymatrix@C6pt@R6pt{
 &           &             &  \c                    &             &     &\\
 &           &             &  \c \ar@{-}[u]            &             &   &\\
\stackrel{\approx}{\De}:&\c\ar@{-}[r]& \c\ar@{-}[r]& \c\ar@{-}[u]\ar@{-}[r]& \c\ar@{-}[r]&\c\ar@{-}[r] &\star\\
}
\end{array}
$$\normalsize
\begin{proposition}
  Let $\La$ be a canonical algebra of domestic representation type $\pP=(p_1,p_2,p_3)$, and $\De$ be the Dynkin diagram $[p_1,p_2,p_3]$. Then the extended canonical algebra $\Lahat=\Lahat\spitz{p_1,p_2,p_3}$ is derived equivalent to the (wild) path algebra  of a quiver $Q$ having extended Dynkin type $\tilde{\tilde{\De}}$.
\end{proposition}
\begin{proof}
Let $Q'$ be a quiver with underlying graph $\tilde{\De}$ such that `the' extension vertex $p$ is a sink. Note that the path-algebra $A'$ of $Q'$ is tilting equivalent to $\La$. Moreover, the one-point extension $A=A'[P_p]$ of $A'$ is isomorphic to the path algebra $A$ of a quiver $Q$ with underlying graph $\tilde{\tilde{\De}}$. Let $\Tt$ be the derived category $\Der{\mmod{A}}$. By construction the `new' indecomposable projective $A$-module $P_*$ corresponding to the one-point extension $A'[P_p]$ is special, and the claim follows.
\end{proof}

\subsubsection*{Euler characteristic zero: the tubular case} Consider a canonical algebra    $\La=\Lapla$ with a tubular weight sequence
$\pP=(p_1,\ldots,p_t)$, we shall assume that $2\le p_1\le p_2\le
\cdots \le p_t$.
\begin{proposition}
The extended canonical algebra
$\Lahat\spitz{\pP,\lala}$ is derived canonical of type $(\bar\pP,\lala)$, where $\bar{\pP}=(p_1,\ldots,p_{t-1},p_t+1)$.
\end{proposition}
\begin{proof}
Let $\XXbar=\XX(\bar{\pP},\lala)$ and $\Tt$ the derived category $\Der{\coh\XXbar}$. Let $E$ be simple in $\coh\XX$ of $\tau$-period $p_t$. Then $E$ is special and $\rperp{E}$ is derived equivalent to $\La(\pP,\lala)$.
\end{proof}
Note that this yields the following wild canonical types $(2,3,7)$, $(2,4,5)$, $(3,3,4)$ and $(2,2,2,3;\la)$. Hence the derived category $\Der{\mmod\Lahat}$ is equivalent to $\Der{\coh{\bar\XX}}$, where $\bar\XX$ is the weighted projective line of type $(\bar\pP,\lala)$.

\subsubsection*{Negative Euler characteristic (the wild case)}
For negative Euler characteristic the derived category of modules
over an extended canonical algebra $\Lahat=\Lahat\spitz{\pP,\lala}$ relates to the study of
a $\ZZ$-graded surface singularity $R=\Rpla$ associated with the weighted projective line $\XX$
associated to $\La$.

For the definition of $\Rpla$ and a table of the singularities with three generators, we refer to Section~\ref{ssect:graded_sing}. The next theorem is a deep result, it is due to the work of several authors and documents the important role, extended canonical algebras have in singularity theory. We will prove it in Sections~\ref{sect:singularities} and~\ref{sect:stable}, where also the relevant definitions are given:

\begin{theorem}
Let $\XX=\Xpla$ be a weighted projective line of negative Euler characteristic, let $\Lahat=\Lahat\spitz{\pP,\lala}$ be the corresponding extended canonical algebra and $R=\Rpla$ be the $\ZZ$-graded singularity attached to $\XX$. Then the derived category $\Der{\mmod\Lahat}$ is triangle-equivalent to each of the following triangulated categories:

 $(i)$ The triangulated category of singularities $\Dsinggrad\ZZ{R}$ of $R$;

 $(ii)$ The stable category $\stabCMgrad\ZZ{R}$ of graded Cohen-Macaulay modules over $R$;

 $(iii)$ The stable category $\stabvect\XX$ of vector bundles on $\XX$.~\hfill\qed
\end{theorem}

\section{Further spectral properties} \label{sect:further_prop}

\subsection{Isospectral algebras} \label{ssect:isospectral}

Let $\coxpol{A}$ be the Coxeter polynomial of a finite dimensional algebra $A$. Its set of roots in the complex plane is denoted by $\Roots{\coxpol{A}}$. The set of roots together with their multiplicities is denoted $\Spec{\coxpol{A}}$, or just $\Spec{A}$, and called the \define{spectrum} of $A$. Two algebras are called \define{isospectral} (or cospectral), if they have the same spectrum, that is, the same Coxeter polynomial. In the same spirit we speak of \emph{isospectral graphs} if their characteristic polynomials are the same. Clearly, derived equivalent algebras are isospectral, but in general isospectral algebras are not derived equivalent, as we are going to illustrate by a couple of examples.

\noindent\emph{Wild hereditary tree algebras which are isospectral but not derived equivalent:}
Consider the tree algebras $A_1$ and $A_2$ given by the displayed quivers:
$$
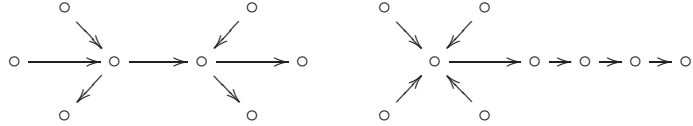
\def\c{\circ}
\xymatrix@R10pt@C8pt{
          & \c\ar[rd]&                & &          &\c\ar[ld]&     \\
  \c\ar[rr]&           &\c\ar[rr]\ar[ld] & &\c\ar[rr]\ar[rd]&         &\c  \\
          &\c         &                & &        &\c       &
}\quad \quad
\xymatrix@R10pt@C8pt{
\c\ar[rd]&         &\c\ar[ld]& &  &     \\
         & \c\ar[rr] &       &\c\ar[r]&\c\ar[r] &\c\ar[r]&\c          \\
\c\ar[ur]&          &\c\ar[lu]&        &       &
}
$$
We denote the corresponding underlying graphs $\De_1$ and $\De_2$. In
\cite{Collatz:Sinogowitz:1957} the graphs $\De_1$ and $\De_2$ were produced as the pair of isospectral graphs with smallest number of vertices, that is, $\charpol{\De_1}(x)=\charpol{\De_2}(x)$. By A'Campo's formula (Proposition~\ref{prop:ACampo}) we have
$\coxpol{A_i}(x^2)= \stcharpol{\De_i}(x)$ for $i=1,2$, hence
$\coxpol{A_1}=\coxpol{A_2}$, that is the algebras $A_1$ and $A_2$ are isospectral. Moreover, we observe that the algebras $A_1$ and $A_2$ are not derived equivalent. Indeed, the quiver $\vDe_i$ appears as a section of a transjective component of the Auslander-Reiten quiver of the derived category of mod $A_i$, for $i=1,2$.

A \define{comb} $[[a_1,a_2,\ldots,a_s]]$ is a tree obtained from a linear `basis', consisting of $s$ consecutive vertices $1,2,\ldots,s$, by attaching to each $i$, $i=1,\ldots,s$, a linear graph $[a_i]$, merging $i$ with an extremal vertex of $[a_i]$. (The attached linear graphs are considered to be disjoint.) The tree $T=[2,2,3,5]$ and the comb $C=[[1,2,2,2,1,1]]$
 $$
 \def\c{\circ}
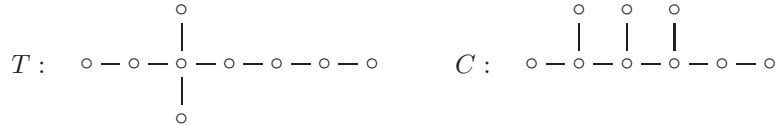
\xymatrix@R9pt@C7pt{ &            &            & \c\ar@{-}[d]           &            &            &            &  &\\ T: &\c\ar@{-}[r]&\c\ar@{-}[r]&\c\ar@{-}[r]&\c\ar@{-}[r]&\c\ar@{-}[r]&\c\ar@{-}[r]&\c&\\ &             &            &\c\ar@{-}[u]           &            &            &            &  &\\
 }
 \quad
 \xymatrix@R9pt@C7pt{
&             &\c\ar@{-}[d]&\c\ar@{-}[d]&\c\ar@{-}[d]&           &   \\
C:& \c\ar@{-}[r]&\c\ar@{-}[r]&\c\ar@{-}[r]&\c\ar@{-}[r]&\c\ar@{-}[r]&\c\\
 }
 $$
are isospectral. By the preceding argument their path algebras $k[\overrightarrow{T}]$ and $k[\overrightarrow{C}]$ will not be derived equivalent, regardless of the chosen orientation.

\medskip\noindent\emph{Isospectral tree algebras with an arbitrary big number of vertices:} Indeed, consider the algebras $A$ and $A'$ given by the following quivers obtained by identifying a vertex of a quiver of type       ${\tilde\EE_8}$ with a extremal vertex of a linear quiver of type $\A{n}$:
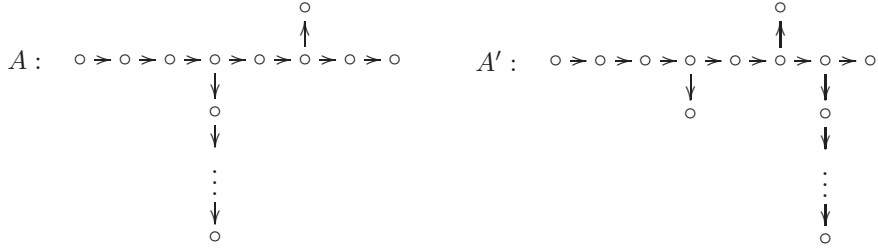
$$
\def\c{\circ}
\xymatrix@R8pt@C6pt{
&          & &    & &                            & \c& &  &  \\
A:&\c\ar[r]& \c\ar[r]&\c\ar[r]&\c\ar[r]\ar[d]&\c\ar[r] &\c\ar[r]\ar[u]& \c\ar[r]&\c\\
&        &         &        &  \c\ar[d]    &         &              &  &        \\        & &         &     &  \vdots\ar[d]&         &              &  &      &  \\
&         &         &     &  \c     &         &              &  &      &  \\
}\quad
\xymatrix@R8pt@C6pt{
&          & &    & &                            & \c& &  &  \\
A':&\c\ar[r]&\c\ar[r]&\c\ar[r]&\c\ar[r]\ar[d]&\c\ar[r]&\c\ar[r]\ar[u]&\c\ar[r]\ar[d]&\c\\
&     &         &      &\c& & & \c\ar[d]      &                                \\
&       &         &        & & & & \vdots\ar[d] &                                \\      &         &         &     & & & & \c          &                                \\
}$$

\begin{lemma}
The algebras $A$ and $A'$ are isospectral.
\end{lemma}
\begin{proof}
It is enough to observe that the underlying graphs $\De$ and $\De'$ satisfy $\charpol\De(x)=\charpol{\De'}(x)$. This follows from the following construction at \cite{Schwenk:1973}:

The \define{coalescence} of $\De_1$ and $\De_2$ at vertices $v_1$ of $\De_1$ and $v_2$ of $\De_2$ is formed by identifying $v_1$ and $v_2$ and denoted by $\De_1 \bullet \De_2$. If $\De_2$ and $\De'_2$ are isospectral graphs and $\De_2\setminus v_2$ and $\De'_2\setminus v'_2$ are also isospectral, then the graphs $\De_1 \bullet \De_2$ and                 $\De_1 \bullet \De'_2$ are isospectral.

To show the claim, only observe that
$$\charpol{\De_1 \bullet \De_2}(x)=\charpol{\De_1}(x)\charpol{\De_2 \setminus v_2}(x)+   \charpol{\De_1\setminus v_1}(x) \charpol{\De_2}(x) -
x \charpol{\De_1\setminus v_1}(x) \charpol{\De_2\setminus v_2}.$$
In our special case $\charpol{\De_2 \setminus v_2}(x)=\charpol{\De'_2\setminus v'_2}(x)=x^2(x^2-2)(x^4-4x^2+2)$.
\end{proof}

\noindent\emph{A Dynkin quiver algebra isospectral to a wild algebra:} The path algebra of a Dynkin quiver of type $\D{12}$ and the extended canonical algebra of type $\spitz{2,4,6}$ are isospectral by Proposition~\ref{prop:anomaly}.

Isospectral problems also illustrate the interplay between spectral graph theory and Coxeter polynomials. The following result
\cite{Lepovic:Gutman:2002}, whose proof we sketch, is an example.

\begin{proposition}
Isospectral stars, with standard orientations, are isomorphic (as quivers or graphs).
\end{proposition}
\begin{proof}
Let $p=[p_1,\ldots,p_t]$ and $q=[q_1,\ldots,q_s]$ be two isospectral stars with the standard orientation. Then $n= 1+\sum_{i=1}^t( p_i-1) = 1+\sum_{j=1}^s( q_j-1)$ and  we assume that $p_1 \ge p_2 \ge \ldots \ge p_t$ and $q_1 \ge q_2 \ge \ldots \ge q_s$. By (1.3), both $\charpol{p}=\charpol{q}$ and $\coxpol{p}=\coxpol{q}$. We shall prove that $p=q$.

First we show that $s=t$. Denote by $c_i$ the coefficient of $x^{n-i}$ in the polynomial $\charpol{p}$.
Since $p$ a tree, the coefficient $c_4$ is the number of pairs of independent edges of $p$, see \cite[Theorem 1.3]{Cvetkovic:Doob:Sachs:1979} . For $p$, an easy computation yields
$$c_4= \begin{pmatrix}n-1 \\2\end{pmatrix} -(n-t-1)-\begin{pmatrix}t \\2 \end{pmatrix}=
\begin{matrix}1 \\2 \end{matrix}[(n-1)(n-4)-t(t-3)].$$
Since, by hypothesis $c_4(p)=c_4(q)$ then $t=s$.

Consider now the expression of the Coxeter polynomial of $p$:
$$\coxpol{p}= [(x+1) -x \sum_{i=1}^t \frac{1-x^{p_i-1}}{1-x^{p_i}}]
\prod_{i=1}^t \frac{1-x^{p_i}}{1-x}.$$
Multiply this polynomial by $(1-x)^t$ to obtain:
$$\Gamma_p=(x+1)\prod_{i=1}^t (1-x^{p_i})- x \sum_{i=1}^t(1-x^{p_i-1})\prod_{j \ne i}(1-x^{p_j})$$
A simple comparison of the coefficients of $\Gamma_p$ and $\Gamma_q$ implies that $p=q$.
\end{proof}

\subsection{Representability of Coxeter polynomials} \label{ssect:representability}
Following \cite{Lenzing:Pena:2006a}, we say that a polynomial $p\in \ZZ[x] $ is
\define{represented} by $q\in \ZZ[x]$ if $p(x^2)=q^*(x):=x^{\deg{q}}q(x+x^{-1})$. It follows that
representable polynomials are self-reciprocal.

The concept of representability arises as a generalization of A'Campo's observation: if $A=k[\vDe]$ is a bipartite hereditary algebra, then $\coxpol{A}$ is represented by   $\charpol{\De}$. We shall see that there are other familiar examples of algebras with representable Coxeter polynomial and will illustrate some applications of this fact.

{\bf Examples:}
Recall  that the (normalized) \define{Chebycheff polynomials} (of
the second kind) $(u_n)_n$ may be inductively constructed by the rules:
$$ u_0=1,\; \ u_1=x, \; \textrm{ and } \;
u_{n+1}=xu_n-u_{n-1} \textrm{ for }n\ge 1.$$
A simple induction shows that the characteristic polynomial of the linear graph
$\AA_n = [n]$ is the  polynomial $u_n$. Moreover
$v_{n+1}$ is represented by $u_n$.

\begin{proposition}
For each $n\ge 2$, the $n$-th cyclotomic polynomial is representable.
In fact, there is an irreducible factor $f$ of $u_{2n-1}$
such that $\Phi_n(x^2)=f^*(x)$.
\end{proposition}

\begin{proof}
For $n=1$ (and $2n-1=1$) we have:
$$u^*_1=x(x+x^{-1})=x^2+1=\Phi_2(x^2).$$
Assume $f^*_i(x)=\Phi_i(x^2)$ for $i<n$, then
$$\Phi_n(x^2)\prod_{\ofrac{1<d<n}{d|n}}\Phi_d(x^2)=v_{n-1}(x^2)=
\frac{x^{2n}-1}{x^2-1}=\frac{v_{2n-1}(x)}{x+1}=
\frac{u^*_{2n-1}(x)}{x+1}.$$
Using (1.4), consider a decomposition
$u_{2n-1}=\prod^s_{i=1}g_i$ in $\ZZ[x]$ such that $g_i$
is irreducible as a self-reciprocal polynomial. Recall that all the cyclotomic polynomials $\Phi_d$ (and hence also $\Phi_n(x^2)$) are irreducible.
By induction hypothesis we may assume that for each $d|n$, we get
$\Phi_d(x^2)=f^*_d(x)$ for certain $f_d=g_i$, $1\le i\le s$.
The result for $\Phi_n$ follows from (1.4).
\end{proof}
As concrete examples we calculate: $\Phi_3=x^2+x+1$ with $\Phi_3(x^2)=u^*_2$; $\Phi_4=x^2+1$
with $\Phi_4(x^2)=f^*_4$ for $f_4=x^2-2$ which is a factor of
$u_3=x^3-2x$ and also of $u_7=x^7-6x^5+10x^3-4x$;
$\Phi_5=x^4+x^3+x^2+x+1$ with $\Phi_5(x^2)=u^*_5(x)$.

The following remark shows the representability of the Coxeter polynomial for many algebras.

\begin{proposition} \label{prop:coxpol_onepoint}
Let $A=B[P]$ be a one-point extension of an algebra $B$ with an
indecomposable projective module $P$ associated to a source $b$ in
$B$. Consider the quotient algebra $C=B/(b)$, then the following holds:

$(a)$  $\coxpol{A}=(1+x)\coxpol{B}-x\coxpol{C}$;

$(b)$ Assume that $\coxpol{B}$ and $\coxpol{C}$ are representable, then $\coxpol{A}$ is representable.
\end{proposition}

\begin{proof}
The formula (a) is shown in Proposition~\ref{prop:onepoint_coxpol}.
Therefore
$$
\coxpol{A}(x^2)=(1+x^2)\coxpol{B}(x^2)-x^2\coxpol{C}(x^2)=x^n(x+x^{-1})q_1(x+x^{-1})-x^n q_2(x+x^{-1})\,,
$$
when $n$ is the degree of $\coxpol{A}$. Hence $\coxpol{A}(x^2)=q^*(x)$ with
$q=xq_1-q_2$.
\end{proof}

\subsubsection*{Graphical representability} In general, we shall say that a polynomial $p\in \ZZ[x]$ is \define{graphically represented} if $p(x^2)=\stcharpol\De(x)$ for a graph $\De$. Applications of graphical representability arise as consequence of the following elementary remark:

\begin{proposition}
Let $p$ and $q\in \ZZ[x]$ be such that $p(x^2)=q^*(x)$. Then

{\rm (a)} If $\mu \in \Roots{p}$, then
$\mu^{-1},\bar{\mu},\bar{\mu}^{-1}\in \Roots{p}$

{\rm (b)} $\Roots{p}\subset \SS^1\cup \RR^+$ (where
$\SS^1=\{v\in \C\colon \Vert v\Vert =1\}$ is the unit circle and
$\RR^+$ the positive real numbers, $0\in \RR^+$) if and only if
$\Roots{q}\subset \RR$.

{\rm (c)} $\Roots{p}\subset \SS^1$ (resp. $\SS^1\setminus
\{1\}$) if and only if $\Roots{q}\subset [-2,2]$ (resp. $(-2,2)$).

{\rm (d)} If $p$ is graphically represented, then
$\Roots{q}\subset \RR$ and $\Roots{p}\subset \SS^1\cup \RR^+$.

{\rm (e)} Assume $p$ is graphically represented by the graph $\De$,
that is $q=\charpol\De (x)$, then

\hspace*{\parindent}
{\rm (i)} $\De$ is a bipartite graph.

\hspace*{\parindent}
{\rm (ii)} $\Roots{p}\subset \SS^1\setminus \{1\}$ if and
only if $\De$ is a union of Dynkin graphs.
\end{proposition}

\begin{proof}
Observe that $0\ne \lambda\in \CC$, then $\mu=\lambda^2\in\Roots{p}$
if and only if $\lambda +\lambda^{-1}\in\Roots{q}$. Assume $\lambda
=r(a+ib)$ with $r\in \RR^+$, $a^2+b^2=1$, then $\lambda +
\lambda^{-1}=(r+r^{-1})a+i(r-r^{-1})b$.

(a): if $\mu =\lambda^2\in\Roots{p}$, then $\mu^{-1}$
yields $\lambda +\lambda^{-1}\in \Roots{q}$, then $\mu^{-1}\in
\Roots{p}$. Since $p$ has real coefficients, then $\bar{\mu}$
(and therefore $\bar{\mu}^{-1}$) is in $\Roots{p}$.

(b): Observe $\mu \in \Roots{p}$ lies in $\SS^1$ (resp. in $\RR^+$)
iff $r=1$ (resp. $b=0$) iff $\lambda +\lambda^{-1}\in \RR$.

(c): In the above case $\lambda +\lambda^{-1}\in [-2,2]$ iff $b=0$
iff $\mu \in \SS^1$.

(d): Assume $q=\charpol\De$ for $\De$ a graph.
Then $\Roots{p}\subset
\SS^1\setminus \{1\}$ iff $\Roots{q}\subset (-2,2)$, in this case
$\De$ is a union of Dynkin graphs, see \cite{Cvetkovic:Doob:Sachs:1979}.
\end{proof}

\begin{Example}
In general the Coxeter polynomial of a path algebra $A=k[\vDe]$ is not graphically representable. Indeed, there are wild quivers $\vDe$ where $\Roots{\coxpol{A}}$ is not contained in $\RR\cup \SS^1$, see~\cite[example~18.1]{Lenzing:cox}. Such a quiver must not be bipartite. Moreover, since $\Roots{\coxpol{A}}$ is closed under complex conjugation and taking inverses, we need at least six vertices for $\vDe$ to make the phenomenon appear. With this information at hand, many examples like
$$
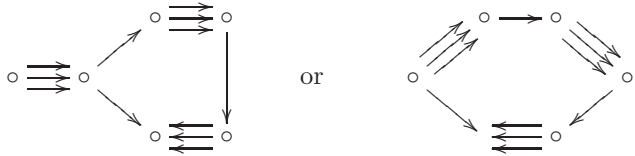
\def\c{\circ}
\xymatrix@R12pt@C16pt{
         &                  &\c\ar[r]\ar@<-1ex>[r]\ar@<1ex>[r]          &\c\ar[dd]                                  &\\
\c\ar[r]\ar@<-1ex>[r]\ar@<1ex>[r]& \c\ar[ru]\ar[rd] &                                           &                                  &\textrm{ or }\\
         &                  &\c                                         &\c\ar[l]\ar@<-1ex>[l]\ar@<1ex>[l] &\\
}
\xymatrix@R12pt@C16pt{
&                                    &\c\ar[r]          &\c\ar[dr]\ar@<-1ex>[dr]\ar@<1ex>[dr]&\\
&\c\ar[ru]\ar@<-1ex>[ru]\ar@<1ex>[ru]\ar[dr]&                  &                                    &\c\ar[dl]                          \\
&                                    &                  &\c\ar[l]\ar@<-1ex>[l]\ar@<1ex>[l]   &\\
}
$$
can be obtained.
\end{Example}

The following result stresses the relationship between properties of the Coxeter transformation $\coxpol{A}$ of an algebra $A$ and the structure of the graph $\De$ in case $\coxpol{A}$ is graphically represented by $\De$.

\begin{proposition}
Let $A$ be an algebra whose Coxeter polynomial $\coxpol{A}$ is graphically represented by $\De$, that is, $\coxpol{A}(x^2) = \stcharpol\De$. Assume that the characteristic polynomial of the graph $\De'$ obtained by deleting a vertex from $\De$ has at most $n - 1$ zeros in the interval $[2,\infty)$. Then $\coxpol{A}$ has at most $2n$ zeros outside the complex unit circle.
\end{proposition}

\begin{proof}
The equation $\coxpol{A}(x^2) = \stcharpol\De(x)$ implies that the cardinality of $\Spec{\coxpol{A}}\cap (\RR^+ \setminus \{ 1 \})$ is twice the cardinality of
$\Spec{\charpol\De} \cap (2,\infty)$. Hence, it is enough to prove that $\charpol\De$ has at most $n$ zeros in the interval
$[2,\infty)$. Indeed, the adjacency matrix $A$ of $\De$ has the shape
$$A= \begin{pmatrix} A'& y^t \cr y & 0 \cr \end{pmatrix},$$
where $A'$ is the adjacency matrix of the graph $\De'$ and $y$ is a row vector.
Consider $\nu_1 \ge \nu_2 \ge \dots \nu _{n-1}$ the roots of $\charpol\De'$ and
$\mu_1 \ge \mu_2 \ge \dots \mu _{n-1} \ge \mu_n$ the roots of $\charpol\De$. The  \define{interlacing property} of eigenvalues yields inequalities:
$$\mu_1 \ge \nu_1 \ge \mu_2 \ge \nu_2 \ge \dots \nu _{n-1} \ge \mu_n$$
which shows the result.
\end{proof}

Several results have been proved which are special cases of the above Proposition. Namely:

(a)~\cite[Proposition 2.6]{Pena:1992} Let $\De$ be a tree graph with $s$ ramification points and $\vDe$ be a quiver with underlying graph $\De$. The Coxeter polynomial $\coxpol{A}$ of the hereditary algebra $A=k[\vDe]$ has at most $2s$ roots outside the unit circle.

(b)~In \cite{Lakatos:2001} \define{generalized stars} are introduced as amalgamations of linear quivers or even cycles in a selected vertex. It is shown that the Coxeter polynomial $\coxpol{S}$ of a wild generalized star $S$ is of the form $p q$ where $p$ is a product of cyclotomic polynomials and $q$ is an irreducible polynomial with exactly one root outside the unit circle (such polynomial is called a \define{Salem polynomial}).

\subsubsection*{Representability of the Coxeter polynomial of an extended canonical algebra} Fix the following notation:
Let $H$ be the star of type $[p_1,\ldots,p_t]$ with Coxeter
polynomial $\coxpol{[p_1,\ldots,p_t]}$ and let $\La$ be a canonical
algebra of weight type $(p_1,\ldots,p_t)$ with
Coxeter polynomial $\coxpol{(p_1,\ldots,p_t)}$; let $\Lahat$ be the
extended canonical algebra with Coxeter polynomial
$\hat{\chi}_{\spitz{p_1,\ldots,p_t}}$.

\begin{proposition}
The above Coxeter polynomials are representable in the following way:

{\rm (a)} $\coxpol{[p_1,\ldots,p_t]}(x^2)=
\charpol{[p_1,\ldots,p_t]}^*$, where $\charpol{[p_1,\ldots,p_t]}$
denotes the characteristic polynomial of the adjacency matrix of star graph
$[p_1,\ldots,p_t]$;

{\rm (b)} $\coxpol{(p_1,\ldots,p_t)}(x^2)=
\stcharpol{K_2}\prod^t_{i=1}\charpol{[p_i-1]}^*$, where
$\charpol{[p]}$ denotes the characteristic polynomial of the linear
path $[p]$ and $K_2$ is the Kronecker diagram $\xymatrix@C16pt@R6pt{\ci
\ar@{-}@/^/[r]&\ci \ar@{-}@/^/[l]}$;

{\rm (c)} $\hat{\chi}_{\spitz{p_1,\ldots,p_t}}(x^2)=
\left(x \charpol{K_2}\prod^t_{i=1} \charpol{[p_i-1]}-
\charpol{[p_1,\ldots,p_t]}\right)^*$.
\end{proposition}

Let $q_{(p_1,\ldots,p_{t})}$ be the polynomial representing
$\hat{\chi}_{(p_1,\ldots,p_t)}$. Under the conditions shown in the above Proposition, a version of \define{Sturm's
Theorem}~\cite{Obreshkoff:1963} assures that given any interval $[\alpha
,\beta]\subset \RR$ and the roots $\lambda_1\le \cdots \le \lambda_s$
of $q_{(p_1,\ldots,p_{t+1})}$ in $[\alpha,\beta]$, then
$q_{(p_1,\ldots,p_t)}$ has roots $\lambda'_1\le \cdots \le
\lambda'_{s-1}$ in $[\alpha,\beta]$ satisfying $\lambda_1\le
\lambda'_1\le \lambda_2\le \lambda'_2\le \cdots \le \lambda_{s-2}\le
\lambda'_{s-1}\le \lambda_s$ (interlacing property). These are the main ingredients for showing the next result~\cite{Lenzing:Pena:2006a}.

\begin{theorem}
Consider the extended canonical algebra $A$ of type
$(p_1,\ldots,p_{t-1},p_t+1)$ and $A'$ an extended canonical algebra
of type $(p_1,\ldots,p_t)$. Then the following holds:

{\rm (a)} If $\specrad{\coxpol{A}}=1$, then also $\specrad{\coxpol{A'}}=1$.

{\rm (b)} $\coxpol{A}$ accepts at most $4$ eigenvalues outside
$\SS^1$.~\hfill\qed
\end{theorem}

\section{One-point extensions} \label{sect:one_point}

\subsection{Fundamental facts} \label{ssect:onepoint_fundamental}
We summarize a number of facts on one-point extensions, starting with two useful general results. The first concerns Hochschild cohomology, the second derived equivalence.

\begin{proposition}[\cite{Happel:1989}] \label{prop:Happel_Hochsch}
If $A=B[M]$ is the one-point extension of $B$ by an module $M$, then the following holds:

$(i)$ There is a long exact sequence $0\ra \hoch{0}{A}\ra \hoch{0}{B}\ra \End{}{M}/k\ra\hoch{1}{A}\ra\hoch{1}{B}\ra\Ext1{}{M}{M}\ra\hoch{2}{A}\ra\cdots $ in Hochschild cohomology.

$(ii)$ If $M$ is exceptional, then $A$ and $B$ have the same Hochschild cohomology in all degrees.~\hfill\qed
\end{proposition}

\begin{proposition}[\cite{Barot:Lenzing:2003}] \label{prop:one-point}
Let $A$ and $B$ be two finite-dimensional $k$-algebras, $M$ an $A$-module and $N$ a $B$-module. Let $\bar{A}$ and $\bar{B}$ be the respective one-point extensions. Then any triangle-equivalence $\Der{\mmod{A}}\ra\Der{\mmod{B}}$ sending $M$ to $N$ yields a triangle-equivalence $\Der\mmod{\bar{A}}\ra\Der\mmod{\bar{B}}$.~\hfill\qed
\end{proposition}
\begin{coro}
Assume $A$ has finite global dimension.
Let $S$ be a simple $A$-module. Let $P(S)$ (resp.\ $I(S)$) denote its projective (resp.\ injective) hull. Then the one-point extensions $A[P(S)]$ and $A[I(S)]$ are derived equivalent.
\end{coro}
\begin{proof}
Let $\tau$ denote the Auslander-Reiten translation in $\Der{\mmod{A}}$. Then the self-equivalence $X\mapsto \tau X[1]$ of $\Der{\mmod{A}}$ maps $P(S)$ to $I(S)$. Here, we view $A$-modules as complexes concentrated in degree zero.
\end{proof}

\subsection{One-point extensions and perpendicular categories} \label{ssect:one_point_perp}
Throughout this subsection we assume that $A$ is a finite dimensional algebra of finite global dimension,
implying that the bounded derived category $\Der{{A}}$ of finite dimensional $A$-modules is homologically finite. Here, a
triangulated category $\Tt$ is called \define{homologically finite} if for any two objects $X$ and $Y$ from $\Tt$ the space $\Hom{\Tt}{X}{Y[n]}$ is non-zero only for finitely many $n$. We will only consider such triangulated categories. Recall that an object $E$ in a triangulated category $\Tt$ is called \define{exceptional} if $\End{}{E}=k$ and, moreover, $E$ has no self-extensions, that is, $\Hom\Tt{E}{E[n]}=0$ for each non-zero integer $n$. Correspondingly, an $A$-module $E$ is called exceptional if $\End{}{E}=k$, and $\Ext{n}{A}{E}{E}=0$ for each integer $n>0$. Note that a module $E$ is exceptional in $\mmod{A}$ if and only if it is exceptional as an object in the triangulated category $\Der{\mmod{A}}$ under the standard embedding from $\mmod{A}$ to $\Der{\mmod{A}}$.

Consider an exceptional object $E$ in a triangulated category, then the \define{right perpendicular category} $\rperp{E}$  of $E$ consists of all objects $X$ from $\Tt$ satisfying the conditions  $\Hom{\Tt}{E}{X[n]}=0$ for each integer $n$. Viewed as a full subcategory of $\Tt$, $\rperp{E}$ is a triangulated category, and the exact inclusion $\rperp{E}\incl \Tt$ admits an exact left adjoint $\ell:\Tt\ra \rperp{E}$, see~\cite{Bondal:Kapranov:1990}.

The next proposition shows that forming one-point extensions is in some sense inverse to forming perpendicular categories (with respect to an exceptional object).
\begin{proposition}
  Assume $A$ has finite global dimension and $E$ is an exceptional object in $\Der{\mmod{A}}$. Then the following holds:

  $(i)$ $\rperp{E}$ has a tilting object $T'$. Hence $\rperp{E}$ is equivalent to $\Der{\mmod{B}}$ for $B=\End{}{T'}$. Moreover, $B$ has finite global dimension.

  $(ii)$ If  $F$ is an exceptional object from $\Der{\mmod{A}}$, satisfying $\Hom{}{T'}{F[n]}=0$ for each $n\neq0$, and $M=\Hom{}{T'}{F}$, then $A$ is derived equivalent to the one-point extension $B[M]$.

$(iii)$ Conversely, assume that $\bar{A}=A[M]$ is the one-point extension of an $A$-module $M$, and let $P$ denote the
indecomposable projective right $\bar{A}$-module $[M,k]$, then $\rperp{P}$ is equivalent to $\Der{\mmod{A}}$.
\end{proposition}
\begin{proof}
  We only sketch the proof. Concerning (i) we note that $T'=\ell(A)$ is a tilting object in $\lperp{E}$, implying the remaining claims. Concerning (ii) the assumptions ensure that $T=T'\oplus F$ is a tilting object for $\Der{\mmod{A}}$, implying that $\End{}{T}=\left[ \begin{array}{cc} B & 0 \\ M & k
  \end{array}\right]$ equals $B[M]$. Now use that $\End{}{T}$ is derived equivalent to $A$, since $T$ is tilting in $\Der{\mmod{A}}$. Assertion $(iii)$ follows, observing that $T'$ is tilting in $\rperp{P}$.
\end{proof}
\begin{remark}
  (i) In a hereditary setting, that is if $\Der{\mmod{A}}$ is equivalent to $\Der\Hh$ for a hereditary category $\Hh$,  the above arguments can be executed within the  category $\Hh$, compare~\cite{Geigle:Lenzing:1991}. For algebras that are not satisfying this condition it is unavoidable to pass to the context of triangulated categories.

  (ii) In practice, it is impossible to describe the algebra $B$ from the proposition in the form $k[\vDe]/I$. The reason is, that the left adjoint $\ell$ does not preserve indecomposability, that is, $B$ is usually not a basic algebra.
\end{remark}
It looks therefore surprising that the above proposition will allow us to calculate the Coxeter polynomial $\coxpol{B}$ of the perpendicular category $\rperp{E}$, as will be shown in the next proposition, which will need some preparation. The \define{Poincar\'{e} series} of an element $e$ in $\Knull{A}$ is defined as $P_{(A,e)}=\sum_{n=0}^\infty \eulerform{e}{\tau^n e}x^n$, where $\tau=\cox{A}$. If $e=[E]$ for some module $E$ we write $P_{(A,E)}$ instead.

\begin{proposition}[\cite{Lenzing:cox}, Prop.~18.3, Cor.~18.2]\label{prop:onepoint_coxpol}
Let $E$ be an exceptional $A$-module and $e$ its class in $\Knull{A}$ and let $B$ be an algebra such that $\mmod{B}$ is derived equivalent to the perpendicular category $\rperp{E}$, formed in $\Der{\mmod{A}}$.  Then the following holds:

$(i)$ The Coxeter polynomial $\coxpol{B}$ is given by the expression
  $
  \coxpol{B}=\coxpol{A}(1-P_{(A,E)})/x.
  $
$(ii)$ The Coxeter polynomial for the one-point extension $\bar{A}=A[E]$ is given by
$\coxpol{\bar{A}}=(1+x)(x+P_{(A,E)})$.

$(iii)$ We have $\coxpol{\bar{A}}=(1+x)\coxpol{A}-x\coxpol{B}$, equivalently $\coxpol{B}=\displaystyle\frac{(1+x)\coxpol{A}-\coxpol{\bar{A}}}{x}$.
\end{proposition}
\begin{proof}
Assertion $(i)$ and $(ii)$ are \cite[Prop.~ 18.3]{Lenzing:cox} and \cite{Lenzing:cox}[Cor.~18.2]. Then $(iii)$ follows from $(i)$ and $(ii)$, eliminating $P_{(A,E)}$.
\end{proof}

\subsection{One-point extensions of hereditary algebras} \label{ssect:one_point:hereditary}
\begin{proposition} \label{cor:dynkin}
Let $A$ be a $k$-algebra which is derived equivalent to the path algebra of a Dynkin quiver $\vDe$. Let $B$ be a one-point-extension or coextension by an indecomposable $A$-module $M$. Then $B$ is derived equivalent to the path algebra of a quiver $\vGa$ obtained from $\vDe$ by adding a new vertex and a new arrow.
\end{proposition}

\begin{proof}
We start with a triangle-equivalence $\varphi$ from $\Der{\mmod{A}}$ to $\Der{\mmod{k[\vDe]}}$. Changing $\varphi$ by a translation $[n]$, if necessary, we may assume that $N:=\varphi(M)$ is a $k[\vDe]$-module. (This uses that $M$ is indecomposable and $A'=k[\vDe]$ is hereditary.) It follows from Proposition~\ref{prop:one-point} that $A[M]$ is derived equivalent to $A'[N]$. Moreover, $N$ is in the $\tau$-orbit (with $\tau$ the Auslander-Reiten translate in $\Der{\mmod{A'}}$) of some indecomposable projective $A'$-module $P(a)$, corresponding to the vertex $a$ of $\vDe$. Invoking the proposition again, we see that $A'[N]$ and $A'[P(a)]$ are derived equivalent. One-point extension with an indecomposable projective $P(a)$ (or injective $I(a)$) turns the path algebra of the quiver $\vDe$ in the path-algebra of the quiver $\vGa$, obtained from $\vDe$ by adding a new arrow (with a new vertex) at $a$.
\end{proof}
Next we deal with the tame situation.
\begin{proposition} \label{prop:one_point_tame_her}
Let $A$ be derived equivalent to the path algebra $H=k[\vDe]$ of an extended Dynkin quiver and consider an exceptional $A$-module $M$. Then the algebra $A[M]$ is derived equivalent to an algebra $B$, which is

$(i)$ wild hereditary if $M$ is derived preprojective or preinjective;

$(ii)$ supercanonical of restricted type if $M$ has derived finite length. Moreover, if $M$ is derived quasi-simple then $B$ is a canonical algebra.
\end{proposition}
\begin{proof}
By Proposition~\ref{prop:one-point} we can assume that $M$ is an exceptional $H$-module which in case (i) is preprojective or preinjective and in  case (ii) is an exceptional regular module from an exceptional tube with a quasi-length strictly less than the rank of the tube. Claim (i) is then proved as in the preceding proposition, and claim (ii) follows from \cite{Lenzing:Pena:2004}, the last one from \cite{Lenzing:Meltzer:1996}.
\end{proof}

\begin{remark}
 Starting with a wild hereditary algebra $H$, a one-point extension $H[M]$ with a module $M$  that is indecomposable preprojective or preinjective is again derived wild hereditary. If $M$ is an exceptional regular $H$-module we don't have a complete picture about the algebras arising as one-point extensions. It is known that there are such one-point extensions which are canonical, or super-canonical or extended canonical, but there are many others, where the `type' is not known. For information on one-point extension algebras of wild hereditary algebras we refer to~\cite{Kerner:Skowronski:2001}, \cite{Kerner:Skowronski:Yamagata:2006} and \cite{Chesne:2004}.
\end{remark}

\subsection{One-point extensions of canonical algebras} \label{ssect:onepoint_canonical}
To study one-point extensions by exceptional modules over canonical algebras, results from \cite{Lenzing:Meltzer:1996}, \cite{Lenzing:Pena:2004},  \cite{Malicki:Skowronski:2005} and \cite{Skowronski:2001} are useful.

Let $\La$ be a canonical algebra of weight type $(p_1,\ldots,p_t)$.

$(a)$ If $M$ is regular simple in the $i$-th exceptional tube $\Tt_i$ (of $\tau$-period $p_i$), then the one-point extension $\La[M]$ is tilting-equivalent to the canonical algebra of weight type $(p_1,\ldots,p_i+1,\ldots,p_t)$ having the same parameter sequence as $\La$.

$(b)$ If $M$ has quasi-length $s$ in $\Tt_i$ (recall this means that $s<p_i$), then $\La[M]$ is derived equivalent to a supercanonical algebra in the sense of \cite{Lenzing:Pena:2004}, where the linear arms of the canonical algebra with index different from $i$ are kept, and the $i$-th arm $1\ra 2 \ra \cdots \ra p_i-1$ is changed to the poset $K(p_i,s)$
\footnotesize$$
\xymatrix@R8pt@C14pt{
       &       &             &                   & \star         &         &    \\
1\ar[r]&2\ar[r]&\cdots\ar[r] &p_i-s\ar[r]\ar[ru] & p_i-s+1\ar[r] &\cdots\ar[r]&p_i-1
}
$$\normalsize
We write $\La(i,s)$ for $\La[M]$ and call it a \define{supercanonical algebra of restricted type}.
This yields the following theorem.
\begin{theorem}\label{thm:1pt-canonical}
Let $A$ be derived equivalent to a canonical algebra $\La$ of weight type $(p_1,\ldots,p_t)$. We fix a triangle-equivalence $\Der{\mmod{A}}=\Der{\coh\XX}$ for the weighted projective line $\XX$ corresponding to $\La$. Let $M$ be an exceptional $A$-module and $\bar{A}=A[M]$. There are two cases:

$(i)$ Up to translation in $\Der{\coh\XX}$ the module $M$ corresponds to an exceptional sheaf of finite length $s$ concentrated in the $i$-th exceptional tube of $\coh\XX$. We have $s< p_i$, and $A[M]$ is derived equivalent to the supercanonical algebra $\La(i,s)$.

$(ii)$ Up to translation in $\Der{\coh\XX}$ the module $M$ corresponds to an exceptional vector bundle over $\XX$. There are three cases:

$(a)$ If $\eulerchar\XX>0$, then $A[M]$ is derived equivalent to the path algebra of a wild connected quiver.

$(b)$ If $\eulerchar\XX=0$, then $(i)$ applies.

$(c)$ If $\eulerchar\XX<0$, and $M$ has rank $\pm1$, then $A[M]$ is derived extended canonical. If $M$ has rank $n$ with $|n|>1$, the derived type of the algebra $A[M]$ is not known in general.
\end{theorem}
\begin{proof} The proof of $(i)$ uses proposition~\ref{prop:one-point} in combination with
\cite{Lenzing:Meltzer:1996} and \cite{Lenzing:Pena:2004}. (ii)(a) is covered by Proposition~\ref{prop:one_point_tame_her}, (ii)(b) follows from \cite{Lenzing:Meltzer:1993}. The first assertion of (ii)(c) is Proposition~\ref{prop:ext_canonical}.
\end{proof}

\section{Accessible algebras} \label{sect:accessible}

In this section we introduce the classes of accessible algebras. Starting from the $k$-algebra $k$, accessible algebras are obtained by successive one-point extensions with exceptional modules. From a homological point of view this class looks quite artificial; in particular, accessible algebras are not closed under Morita equivalence or, more generally, under derived equivalence. It is thus natural to investigate \emph{derived accessible algebras}, that is, those algebras which are derived equivalent to an accessible algebra.

There are many reasons to investigate such algebras: $(i)$ By design accessible algebras have a combinatorial flavor, inherited from the largely combinatorial character of exceptional modules; $(ii)$ Many interesting algebras are derived equivalent to accessible algebras, as will be shown in this section; $(iii)$ Accessible algebras have a strong affinity to spectral analysis; $(iv)$ Derived accessible algebras seem to be those finite dimensional algebras with the closest connection to singularity theory; we refer to Section~\ref{sect:stable} for further details.

\subsection{Basic properties} \label{ssect:basic_prop}

A finite dimensional $k$-algebra $A$ is called \define{accessible} (resp.\ \define{weakly accessible}) if there are finite dimensional algebras $A_1,\ldots,A_n$ such that $A_1=k$, $A_n=A$ and moreover $A_{s+1}$ is (isomorphic to) a one-point extension or coextension of $A_{s}$ with an exceptional (resp.\ an indecomposable) module  $E_{s}$. Weakly accessible algebras (under the name \emph{weakly simply connected}) have been studied in~\cite{Geiss:2002}. One reason to focus on accessible  instead of weakly accessible algebras is the combinatorial flavor of accessible algebra, supported by the following result.

\begin{proposition}
The Hochschild cohomology of a derived accessible algebra $A$ satisfies $\hoch{0}{A}=k$ and $\hoch{i}{A}=0$ for $i\geq1$. In particular, the degree one coefficient of the Coxeter polynomial $\coxpol{A}$ equals one.
\end{proposition}
\begin{proof}
The first assertion follows from Happel's exact sequence, see Proposition~\ref{prop:Happel_Hochsch}, implying that an accessible algebra $A$ has the Hochschild cohomology of $k$. The last assertion then follows from Happel's trace formula (Proposition~\ref{prop:Happel_trace}).
\end{proof}
\begin{Examples} \label{ex:accessible}
  $(i)$ Each path algebra of a tree is accessible.

  $(ii)$ Not every hereditary algebra is accessible. Obviously the Kronecker algebra $k[\ci\rightrightarrows\ci]$ is not accessible. More generally, any algebra $A$ derived equivalent to a hereditary algebra $H$ of extended Dynkin type $\widetilde{\AA}_{p,q}$, equivalently to a canonical algebra of type $(p,q)$, is not accessible. Indeed, the Coxeter polynomial $\coxpol{A}$ has the form $(x-1)^2(1+x+\cdots+x^{p-1})(1+x+\cdots+x^{q-1})$. Hence the degree one coefficient of $\coxpol{A}$ equals zero which is impossible for a derived accessible algebra.

  $(iii)$ Not every poset algebra is accessible: By $(ii)$ the $n$-crown $C_n$, given by the quiver
  $$
 \def\c{\circ}
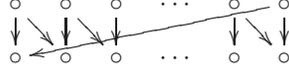
\xymatrix@R10pt@C8pt{
\c\ar[d]\ar[rd]&\c\ar[d]\ar[rd]&\c\ar[d] &\cdots       &\c\ar[d]\ar[rd]&\c\ar[d]\ar[dlllll]\\
\c             &\c             &\c       &\cdots       &\c              &\c
}
  $$
  with $2n$ vertices, is not accessible for $n\geq2$, since it is derived equivalent to the canonical algebra of type $(n,n)$.

  $(iv)$ The path algebra $A$ of the quiver $\vDe$
  $
  \def\c{\circ}
  \xymatrix@R8pt@C6pt{
                  & \c\ar[ld]\ar[rd] &    \\
  \c\ar[d]\ar[drr] &                  &\c\ar[d]\ar[lld]\\
  \c              &                  & \c
  }
  $
  is not accessible. Equipped with all commutativity relations, $\vDe$ yields a poset $P$ whose poset algebra $P$ is accessible

  $(v)$ A canonical algebra $A$ with \emph{at least four weights} is not derived accessible. Indeed, $A$ has the Coxeter polynomial $\coxpol{A}=(x-1)^2\prod_{i=1}^t(1+x+\cdots+x^{p_i-1})$, hence the degree one coefficient of $\coxpol{A}$ equals $t-2\geq2$, which is impossible for an accessible algebra.
\end{Examples}
The above examples provide the proper framework for our next result.
\begin{proposition}
  Each canonical algebra $A$ with three weights is accessible.
\end{proposition}

\begin{proof}
  Notice that $A$ is a one-point extension $A=B[M]$ of the path algebra $B$ of a star $[p_1,p_2,p_3]$ (with each $p_i\geq2$) by an indecomposable module $M$ of dimension vector $[M]=$\scriptsize$\begin{array}{rl}
     &11\cdots 1\\
   2 &11\cdots 1\\
     &11\cdots1
  \end{array}$.\normalsize As is easily checked,  we have $\eulerform{[M]}{[M]}_B=1$ and $\End{}{M}=k$. Since $B$ is hereditary, this implies that $\Ext1{B}{M}{M}=0$. Since $B$ is accessible, the claim follows.
\end{proof}

\subsection{Tree algebras} \label{ssect:tree_algebras}
We start with a lemma that later will be generalized.
\begin{lemma}\label{lemma:linear_quiver}
Each $k$-algebra $A$ given by a  linear quiver
\small$$
\xymatrix{
1\ar[r]^{x_1}  &  2\ar[r]^{x_2}  &\cdots\ar[r]^{x_{n-2}} &n-1\ar[r]^{x_{n-1}}&n
}
$$\normalsize
with zero relations is accessible.
\end{lemma}

\begin{proof} For $s$ from 1 to $n$ let $A_s$ be the restriction of $A$ to the subquiver consisting of the vertices $1,\ldots,s$. By recursion it is sufficient to show that $A$ is isomorphic to the one-point extension of $A_{n-1}$ by an exceptional $A_{n-1}$-module $E$. There are two cases to consider. If no zero relation ends in vertex $n$ then let $E$ be the indecomposable projective $A_{n-1}$-module $P(n-1)$ corresponding to the vertex $n-1$. Otherwise --- assuming an irredundant set of zero relations of length $\geq2$ --- we have a zero relation $\xymatrix{r\ar[r]^{x_{n-1}\cdots x_r}&n}$ for a uniquely defined $r<n-1$, and let $E$ be the indecomposable injective $A_{n-1}$-module $I(r+1)$ corresponding to the vertex $r+1$.
\end{proof}

Since the algebras given by a linear quiver with zero relations are obviously representation-finite they may not look to be very interesting. We will show, however, that the closure of this class against derived equivalence will contain many interesting algebras mostly of wild representation type. Often this insight is due to spectral analysis. Assume for this, more specifically, that $A_n$ denotes the linear quiver
\small$$
\xymatrix{
1\ar[r]^x  &  2\ar[r]^x  &\cdots\ar[r]^x &n-1\ar[r]^x&n
}
$$\normalsize
satisfying all zero relations $x^3=0$. An easy calculation shows that $A_{11}$ and the wild canonical algebra $\La$ of weight type $(2,3,7)$ are isospectral (with Coxeter polynomial $(x-1)^2v_2v_3v_7$ and the same holds true for $A_{12}$ and the extended canonical algebra $\Lahat\spitz{2,3,7}$ which both have Coxeter polynomial $\Phi_{42}$. This poses the question whether the algebras in question are derived equivalent. In fact they are, as we are going to show later in this section.

Path algebras of linear quivers with zero relations are a special case of so called tree algebras. By definition a \define{tree algebra} arises from the path algebra of a (finite) tree by factoring out zero relations.

\begin{proposition}
  Each tree algebra is accessible.
\end{proposition}
 \begin{proof}
Let $A=k[\vDe]/I$ be a tree algebra and choose a terminal vertex $a$ such that the quotient $B=A/(a)$ is connected. We may assume that $A=B[M]$ for an indecomposable $B$-module $M$. Clearly, the support algebra of $M$, formed by those vertices $x$ in $\vDe$ where $M(x) \ne 0$, is a connected hereditary algebra
$C= k [\vDe']$ which is convex in $A$, that is $\vDe'$ is path closed in $\vDe$. The result follows by induction using the following two claims:

(i) Let $M$ be the unique indecomposable module over a path tree algebra $C=k[\De]$ with dim$_kM(x)=1$ for every vertex $x$. Then $M$ is exceptional.

(ii) Let $C$ be a convex subcategory of a triangular algebra $A$ and consider two $C$-modules $X, Y$. Then Ext$_C^i(X,Y)=$ Ext$_A^i(X,Y)$ for any $i \ge 0$.

For (i) observe that, since $\De$ is a tree $q_C([M])=1$ where
$$q_C(v)= \sum_{i \text{ vertex}}v(i)^2 - \sum_{i \ra j} v(i)v(j)$$ is the \define{Tits form} of $C$ which, due to heredity of $C$ agrees with the Euler quadratic form. Moreover $q_C([M])=$ $\ddim{k}{\End{C}{M}}-\ddim{k}{\Ext1{C}{M}{M}}$ and $\End{C}{M}$ is trivial. Hence $M$ is exceptional. For (ii), to proceed by induction, we may assume that $A=C[M]$, then for $i \ge 1$ we have that
$\Ext{i}{C}{X}{Y}=\Ext{i-1}{C}{K}{Y}=\Ext{i-1}{A}{K}{Y}$, where $0 \ra K \ra P \ra X \ra 0$ is an exact sequence with $P$ a projective $C$-module. Since $P$ is also
$A$-projective, we get the result.
 \end{proof}

\subsection{Poset algebras} \label{ssect:poset_algebras}
Many \define{poset algebras}, that is, path algebras of finite fully commutative quivers, are accessible. Adapting the argument from Lemma~\ref{lemma:linear_quiver}, it is easily seen that each poset of type
\begin{equation} \label{eq:ladder1}
\def\ci{\circ}
B_n\quad \xymatrix@R12pt@C10pt{
 \ci\ar@{.}[dr] \ar[d]\ar[r] &\ci\ar@{.}[dr] \ar[d]\ar[r] & \ci\ar[d]\ar[r] &\cdots\ar[r]&\ci\ar@{.}[dr]\ar[d]\ar[r]&\ci\ar[d]\\
\ci \ar[r]       &\ci \ar[r]       &\ci\ar[r] &\cdots\ar[r]& \ci\ar[r]&\ci\\}
\end{equation}
leads to an accessible algebra. This example allows variations like
\begin{equation} \label{eq:ladder2}
\def\ci{\circ}
C_n\quad \xymatrix@R12pt@C10pt{
 \ci\ar@{.}[dr] \ar[d]\ar[r] &\ci\ar@{.}[dr] \ar[d]\ar[r] & \ci\ar[d]\ar[r] &\cdots\ar[r]&\ci\ar@{.}[dr]\ar[d]\ar[r]&\ci\ar[d]\ar[r]&\ci\\
\ci \ar[r]       &\ci \ar[r]       &\ci\ar[r] &\cdots\ar[r]& \ci\ar[r]&\ci\\}
\end{equation}
or
\begin{equation} \label{eq:ladder3}
\def\ci{\circ}
D_n\quad\xymatrix@R12pt@C10pt{
 \ci\ar@{.}[dr] \ar[d]\ar[r] &\ci\ar@{.}[dr] \ar[d]\ar[r] & \ci\ar[d]\ar[r] &\cdots\ar[r]&\ci\ar@{.}[dr]\ar[d]\ar[r]&\ci\ar[d]\ar@{.}[dr]\\
\ci \ar[r]       &\ci \ar[r]       &\ci\ar[r] &\cdots\ar[r]&\ci\ar[r]& \ci\ar[r]&\ci\\}
\end{equation}
mixing commutativity and zero relations. In each case the subscript $n$ denotes the number of vertices. These diagrams (not the algebras!) play a prominent role in singularity theory, see~\cite{Gabrielov:1973} for their first appearance, and \cite{Ebeling:2007} for a recent textbook account. We note that Ladkani~\cite{Ladkani:2007}, \cite{Ladkani:2008} has investigated poset algebras recently.

By way of example we note that the algebras $A_{10}$, $B_{11}$ and $A_{12}$ have Coxeter polynomials $(x-1)^2v_2v_3v_6$, $(x-1)^2v_2v_3v_7$ and $\Phi_{42}$ indicating a possible derived equivalence to canonical $(2,3,6)$, $(2,3,7)$ and extended canonical $\spitz{2,3,7}$, respectively. Again, this is not coincidental and the derived equivalence actually holds.

Alternatively, accessibility of the poset algebras, discussed above, can be derived from our next result.

\begin{proposition}
  For a poset algebra $A$ of a finite poset $S$ the following are equivalent:

  $(i)$ $S$ does not contain convexly any $n$-crown for $n\ge 2$.

  $(ii)$ The poset algebra $A'$ of each convex subset $S'$ of $S$ has $\hoch{1}{A'}=0$.

Moreover, these conditions imply that $A$ is accessible.
\end{proposition}
\begin{proof}
The equivalence of (i) and (ii) is proved in
\cite{Draexler:1994}.

Assume that $A$ satisfies (i) and (ii) and consider $A$ as a one-point extension $A=B[M]$ of a connected algebra $B$. By induction hypothesis, $B$ is accessible. We shall show that $M$ is exceptional. The paper
\cite{GPPRT:1999} shows that all higher Hochschild groups satisfy
$\hoch{i}{A'}=0=\hoch{i}{A}$. Happel's long exact sequence, see Proposition~\ref{prop:Happel_Hochsch} implies the result.
\end{proof}
As shown by Example~\ref{ex:accessible} the conditions $(i)$ and $(ii)$ are only sufficient conditions and don't give a full characterization of accessible poset algebras.

Following \cite{Skowronski:1996a} an algebra $A$ satisfying that each convex subcategory $C$ has vanishing first Hochschild cohomology group is called a \define{strongly simply connected algebra}. From subsection~\ref{ssect:onepoint_fundamental} follows that strongly simply connected algebras are weakly accessible.
For the next remark recall that a \define{schurian algebra} $A=k[\vDe]/I$ satisfies
that $\ddim{k}{A(x,y)} \le 1$ for every couple of vertices $x,y$. In particular, tree algebras and poset algebras are schurian.

\begin{proposition}
Let $A$ be a schurian algebra. The following are equivalent:

\emph{(a)} $A$ is strongly simply connected;

\emph{(b)} for every convex subcategory $C$ of $A$, the algebra $C$ is accessible;

\emph{(c)} for every convex subcategory $C$ of $A$, the Hochschild cohomology groups
$H^i(C)$ vanish, for all $i \ge 1$.
\end{proposition}
\begin{proof}
The equivalence of (a) and (c) was shown in \cite{GPPRT:1999}, while we have shown already that (b) implies (c).
To show that (c) implies (b) it is enough to consider a one-point extension $A=B[M]$ where $A$ satisfies (c) and $B$ is
accessible and show that $M$ is exceptional. It follows that the Hochschild cohomology groups of $B$ and $A$ of degree $\ge 1$ vanish.
Then $M$ is exceptional by an obvious consideration of the long exact Happel's sequence.
\end{proof}

\subsection{Spectral analysis} \label{ssect:spectral_analysis}
We illustrate the methods by investigating the linear quiver with 12 vertices subject to all zero-relations $x^3=0$.
\begin{proposition}
Let $A_n$ be the linear quiver of $n$ vertices with the relations $x^3=0$. Then $A_{11}$ is derived equivalent to the canonical algebra $\La(2,3,7)$ and $A_{12}$ is derived equivalent to the extended canonical algebra $\Lahat\spitz{2,3,7}$.
\end{proposition}

\begin{proof}
The following table displays the Coxeter polynomial of $A_n$ and the corresponding derived type for $n=1,\ldots,12$:
$$
\begin {array}{|l|c|c|c|c|c|c|c|}\hline n&1&2&3&4&5&6&7\\ \hline  
\coxpol{A_n}&v_{{2}}&v_{{3}}&v_{{4}}&{\frac {v_{{2}}v_{{6}}}{v_{{3}}}}&{\frac {v_{{2}}v_{
{8}}}{v_{{4}}}}&{\frac {v_{{3}}v_{{12}}v_{{2}}}{v_{{4}}v_{{6}}}}&{
\frac {v_{{2}}v_{{18}}v_{{3}}}{v_{{6}}v_{{9}}}}\\\hline 
\textrm{type}&[1]&[2]&[2,2]&[2,2,2]
&[2,2,3]&[2,3,3]&[2,3,4]\\\hline
\end {array}
$$
$$
\begin {array}{|l|c|c|c|c|c|}\hline n&8&9&10&11&12\\\hline \coxpol{A_n}&  \Phi_{30}&\left (x-1
\right )^{2}v_{{2}}v_{{3}}v_{{5}}&\left (x-1\right )^{2}v_{{2}}v_{{3}}
v_{{6}}&\left (x-1\right )^{2}v_{{2}}v_{{3}}v_{{7}}&\Phi_{{42}}
\\\hline
\textrm{type}&[2,3,5]&[2,3,6]=(2,3,5)&(2,3,6)&(2,3,7)&\spitz{2,3,7}\\\hline
\end {array}
$$
The table should be read as follows. First the Coxeter polynomials $\coxpol{A_n}$ are determined. Then, inductively, the derived shape of the algebra follows using the preceding analysis of the shape of one-point extensions by exceptional modules. The derived shape for the algebras $A_1,\ldots,A_8$ follows inductively from Corollary~\ref{cor:dynkin}, observing that a one-point extension of a Dynkin diagram is either a Dynkin diagram or an extended Dynkin diagram, the type of the extension is thus determined by its Coxeter polynomial. This also works to determine the type of $A_9$. Here we had two equivalent choices: the extended Dynkin type $[2,3,6]$ and the canonical type $(2,3,5)$. For the following extensions it is preferable to deal with the canonical type. The extension steps from $A_9$ to $A_{10}$ and then to $A_{11}$ are covered by Theorem~\ref{thm:1pt-canonical}, parts $(i)$ and $(ii)$. The extension step from $A_{11}$ to $A_{12}$ is more difficult, since $A_{11}$ is derived equivalent to the wild canonical algebra of type $(2,3,7)$, where we have insufficient knowledge on the one-point extensions, even with exceptional modules. Here we use that the exceptional $A_{11}$-module $M$ with support $\{11,12\}$ has rank one, yielding an extension $A_{12}=A_{11}[M]$ which is derived equivalent to the extended canonical algebra of type $\spitz{2,3,7}$. Note that --- up to a factor minus one --- the rank is given in terms of the Euler form by $\eulerform{-}{w}$, where $w$ is a generator of the radical of the quadratic form.
\end{proof}
The same proof works for other situations as well. By way of example let $A_n$ be the algebra $B_n$ of (\ref{eq:ladder1}) if $n$ is even and equal to the algebra $C_n$ of (\ref{eq:ladder2}) if $n$ is odd. We obtain exactly the same table (and rank information for the extension step from $A_{11}$ to $A_{12}$) as in the above proof and hence the following result:
\begin{coro} The algebra $C_{11}$ of (\ref{eq:ladder2}) is derived canonical of type $(2,3,7)$ and the algebra $B_{12}$ of (\ref{eq:ladder1}) is derived extended canonical of type $\spitz{2,3,7}$ ~\hfill\qed\end{coro}

The preceding discussion shows that the class of accessible algebras is very suitable to spectral analysis. One might therefore ask whether an accessible algebra is determined --- up to derived equivalence --- by its Coxeter polynomial. An interesting counterexample is the following:

Investigating the Coxeter polynomials of linear quivers with zero relations, an interesting anomaly happens for $12$ vertices. Allowing only zero relations involving at least three (consecutive) arrows, the number of such algebras is given by the $n$-th Catalan number $c_n=1/n \binom{2n-2}{n-1}$. For $n=12$ we thus have $8524$ such algebras, yielding $176$ different Coxeter polynomials, one of them the Coxeter polynomial $v_2v_{22}/v_{11}=\Phi_2^2\Phi_{22}$ of the Dynkin quiver $\D{12}$. The number $\de_n$ of such algebras with a Coxeter polynomial of type $\DD$ is given as follows
$$
\begin{array}{|c|ccccc|}\hline
n & 10 & 11 & 12 & 13 & 14 \\\hline
\de_n &7 & 6& 737 & 7 & 7 \\\hline
\end{array}
$$
Focussing on the case of 12 vertices, there are only 6 algebras of the 737-list which are derived equivalent to $\D{12}$, namely four $\{[2, 5], [3, 6] \}$, $\{[3, 6], [4, 7] \}$, $\{[4, 7], [5, 8]\}$, $\{[1, 4], [2, 5]\}$ given by two zero relations and two $\{[1, 4]\}$, $\{[1, 12]\}$ given each by a single relation. The remaining 731 cases may be analyzed by spectral analysis. For instance, the set of four zero relations $\{[1, 6], [3, 8], [6, 9], [7, 12]\}$ yields an extended canonical algebra of type $\spitz{2,4,6}$:
$$
 \begin {array}{|c|c|c|c|c|c|c|c|c|}\hline n& 1&2&3&4&5&6&7&8\\\hline
 \coxpol{A_n}&v_{{2}}&v_{{3}}&v_{{4}}&v_{{5}}&v_{{6}}&\dis{\frac {v_{{2}}v_{{10}}}{v_{{5}}}}&\dis{
\frac {v_{{2}}v_{{18}}v_{{3}}}{v_{{6}}v_{{9}}}}  &\dis\left (x-1\right )^{2}v_{{2}}v_{{3}}v_{{4}}\\\hline
\textrm{type}& [1]&[2]&[3]&[4]&[5]&[2,2,4]&[2,3,4]& (2,3,4)\\\hline
\end {array}
$$
$$
 \begin {array}{|c|c|c|c|c|}\hline n&9&10&11&12\\\hline
\coxpol{A_n}&\dis \left (x-1\right )^{2}v_{{2}}{v_{{4}}}^{2}&\dis\left (x-1\right )^{2}v_{{2}}v_{{4}}v_{{5}}&\dis\left (x-1\right )
^{2}v_{{2}}v_{{4}}v_{{6}}&\dis{\frac {v_{{2}}v_{{22}}}{v_{{11}}}}={\Phi_{{2}}}^{2}\Phi_{{22}}\\\hline
\textrm{type}&(2,4,4)&(2,4,5)&(2,4,6)&\spitz{2,4,6}\\\hline\end {array}
$$
We have shown the next proposition.
\begin{proposition}\label{prop:anomaly}
  The path algebra of a Dynkin quiver of type $\D{12}$ and the extended canonical algebra $\La\spitz{2,4,6}$ are isospectral.~\hfill\qed
\end{proposition}

\section{Singularities} \label{sect:singularities}

This section and the next one deal with one of the lucky instances in mathematics where one has very different descriptions for the same mathematical object, allowing to merge the knowledge from the various perspectives. In our situation we will have --- at least --- four different descriptions of our object of study, all attached to a weighted projective line $\XX$. Three of these objects will be discussed in this section: the triangulated category of singularities $\Dsinggrad\ZZ{R}$ of a graded Gorenstein algebra $R$ attached to $\XX$,  the stable category $\stabCMgrad\ZZ{R}$ of graded Cohen-Macaulay modules over $R$, and the derived category of finite dimensional modules over a finite dimensional algebra $A$, which --- depending on the Euler characteristic of $\XX$ --- is either the path algebra of a Dynkin quiver, or a canonical algebra of tubular type, or, finally, an extended canonical algebra. The fourth object is the stable category $\stabvect\XX$ of vector bundles on $\XX$ that will be discussed in the next section. Orlov's theorem shows that two of the four objects are equivalent, the remaining equivalences are due to the work of several people. To discuss the setup properly, we first review some properties of weighted projective lines and their associated graded singularities, then we discuss Orlov's theorem and finally we apply it to weighted projective lines.

\subsection{The graded singularities associated to a weighted projective line} \label{ssect:graded_sing}

For a given weight sequence $\pP=(p_1,\ldots,p_t)$, $t\geq0$, of integers $p_i\geq2$ we form the rank one abelian group $\LL=\Lp$ on generators $\vx_1,\ldots,\vx_t$ subject to the relations $\vc:=p_1\vx_1=\cdots=p_t\vx_t$. If additionally $\lala=(\la_3,\ldots,\la_t)$ is a parameter sequence of pairwise distinct non-zero elements from $k$, we form the commutative affine $k$-algebra $S=\Spla$ on generators $x_1,\ldots,x_t$ subject to the $(t-2)$ relations $x_i^{p_i}=x_2^{p_2}-\la_ix_1^{p_1}$, $i=3,\ldots,t$. Attaching $x_i$ degree $\vx_i$, the algebra $S$ becomes an $\LL$-graded algebra whose homogeneous components $S_{\vx}$, $\vx\in\LL$,  are finite dimensional $k$-vector spaces.  The \define{weighted projective line} $\XX=\XX(\pP,\lala)$, introduced in \cite{Geigle:Lenzing:1987}, has a category of coherent sheaves $\coh\XX$ which is equivalent to the quotient category
$
\modgrad\LL{S}/\modgradnull\LL{S}
$
of the category $\modgrad\LL{S}$ of finitely generated $\LL$-graded $S$-modules modulo its Serre subcategory $\modgradnull\LL{S}$ of $\LL$-graded $S$-modules of finite length, see \cite{Geigle:Lenzing:1987,Geigle:Lenzing:1991}. If $\Oo$ denotes the structure sheaf of $\XX$, then naturally $\Hom{}{\Oo(\vx)}{\Oo(\vy)}=S_{\vy-\vx}$.

If the (orbifold) \emph{Euler characteristic} $\eulerchar{\XX}=2-\sum_{i=1}^t\left(1-1/p_i\right)$ of $\XX$ is different from zero, then there exists an alternative descriptions of $\coh\XX$ as the quotient category $\modgrad\ZZ{R}/\modgradnull\ZZ{R}$ by a positively $\ZZ$-graded Gorenstein algebra $R$ as follows. Restricting the grading of $S$ to the subgroup of $\Lp$ generated by the dualizing element $\vom=(t-2)\vc-\sum_{i=1}\vx_i$ yields a positively $\ZZ$-graded algebra $R=\Rpla$, $R=\Dir_{n=0}^\infty R_n$, where
$$
R_n=
\begin{cases} S_{-n\vom}\quad \textrm{if}\ &\eulerchar\XX>0\\
S_{n\vom} \quad \textrm{if}          & \eulerchar\XX<0.
\end{cases}
$$
We also need some results, relating the category $\vect\XX$ of vector bundles over $\XX$ to suitable categories of graded maximal Cohen-Macaulay modules over $R$. Recall that a commutative $H$-graded algebra $R$ is called \define{Gorenstein} if the $R$-module $R$ has finite injective dimension. By $\CMgrad{H}{R}$ we denote the category of maximal \define{Cohen-Macaulay modules}. The next statement combines several results: \cite[thm.~5.1]{Geigle:Lenzing:1987}, \cite[ prop.~8.4]{Geigle:Lenzing:1991} and \cite[cor. 5.6 and 5.8]{Lenzing:1994}.

\begin{theorem} \label{thm:gorenstein}
Let $\XX=\Xpla$ be a weighted projective line of Euler characteristic $\eulerchar\XX=2-\sum_{i=1}^t\left(1-1/p_i\right)$. Then the following holds:

$(i)$ For arbitrary Euler characteristic, the $\LL$-graded algebra $S$ is complete intersection, in particular graded Gorenstein.

Moreover, sheafification $M\mapsto \wtilde{M}$ induces a natural equivalence, with inverse the graded global sections functor $\Ga^*=\Dir_{\vx\in\Lp}\Hom{}{\Oo(-\vx)}{-}$, between the categories of $\CMgrad{\LL}{S}$ of graded maximal Cohen-Macaulay modules over $S$ and the category $\vect{\XX}$ of vector bundles on $\XX$.

$(ii)$ If $\eulerchar\XX>0$ then $R$ is a positively $\ZZ$-graded algebra with three generators and a single relation. Accordingly, $R$ is complete intersection, in particular Gorenstein.

Moreover, sheafification $M\mapsto \wtilde{M}$ induces a natural equivalence  $\CMgrad{\ZZ}{R}\ra \vect\XX$ with inverse the global section functor $\Ga^*=\Dir_{n\in\ZZ}\Hom{}{\Oo(n\vom)}{-}$.

$(iii)$ If $\eulerchar\XX<0$ then $R$ is $\ZZ$-graded Gorenstein, in general not complete intersection.

Moreover, sheafification $M\mapsto \wtilde{M}$ induces a natural equivalence  $\CMgrad{\ZZ}{R}\ra \vect\XX$ with inverse the global section functor $\Ga^*=\Dir_{n\in\ZZ}\Hom{}{\Oo(-n\vom)}{-}$.~\hfill\qed
\end{theorem}

Note that $\eulerchar\XX>0$ holds if and only if $t\leq3$ and moreover, the star $[p_1,\ldots,p_t]$ is Dynkin, that is, we deal with one of the weight types $(\;)$, $(p)$, $(p,q)$, $(2,2,n)$, $(2,3,3)$, $(2,3,4)$ or $(2,3,4)$. In this case $R=R(\pla)$ only depends on the weight sequence, the choice of the parameters does not matter. The algebra $R$ has the form $R=k[x,y,z]/(f)$, where the generators $x,y,z$, the relation $f$ and the degrees of generators and relation are given by the table below, see \cite[prop. 8.4]{Geigle:Lenzing:1991}.
\begin{center}\small
\renewcommand\arraystretch{1.4}
\begin{tabular}{|p{60pt}|p{80pt}|p{70pt}|p{90pt}|}\hline
 \textrm{weight type}  & \textrm{generators} $(x, y, z)$  &
\textrm{relation}& \textrm{degrees of} $x,y,z;f$\\
\hline
 $(p, q)$  & $({x}_{0}\,{x}_{1} , {x}_{1}^{p + q} , {x}_{0}^{p + q} )$& ${x}^{p + q}  - y\, z$&
$(1, p, q; p+q)$    \\
\hline
 $(2, 2, 2l)$  & $({x}_{2}^{2} , {x}_{0}^{2} , {x}_{0}\,{x}_{1}\,{x}_{2} )$& ${z}^{2}  + x({y}^{2}  + y\,{x}^{l} )$
& $(2, l, l + 1; 2(l+1))$  \\
\hline
 $(2, 2, 2l + 1)$  & $({x}_{2}^{2} , {x}_{0}\,{x}_{1} , {x}_{0}^{2}\,{x}_{2}
)$  & ${z}^{2}  + x({y}^{2}  + z\,{x}^{l} )$& $(2, 2l + 1, 2l + 2; 4(l+1))$  \\
\hline
 $(2, 3, 3)$  & $({x}_{0} , {x}_{1}\,{x}_{2} , {x}_{1}^{3} )$  & ${z}^{2}  + {y}^{3}  + {x}^{2}\, z$& $(3, 4, 6;12)$
 \\
\hline
 $(2, 3, 4)$  & $({x}_{1} , {x}_{2}^{2} , {x}_{0}\,{x}_{2} )$& ${z}^{2}  + {y}^{3}  + {x}^{3}\, y$  & $(4, 6, 9; 18)$
 \\
\hline
 $(2, 3, 5)$  & $({x}_{2} , {x}_{1} , {x}_{0} )$&
${z}^{2}  + {y}^{3}  + {x}^{5}$  & $(6, 10, 15; 30)$  \\
\hline
\end{tabular}
\end{center}
We note that the Poincar\'{e} series $P_R=\sum_{n=0}^\infty \ddim{k}{R_n}x^n$ is given by the formula
$$
P_R=\frac{1}{(1-x)^2}\frac{\coxpol{[p_1,p_2,p_3]}}{\coxpol{(p_1,p_2,p_3)}},
$$
and a similar remark applies to the table below.
If $\eulerchar\XX=0$, then $\XX$ is called a \define{tubular curve} and the weight type is --- up to ordering --- one of $(2,3,6)$, $(2,4,4)$, $(3,3,3)$ or $(2,2,2,2)$. In the tubular case there does not exist a  $\ZZ$-graded Gorenstein algebra serving as a graded coordinate algebra for $\coh\XX$.

In the wild case there are exactly 14 cases, where $R$ is generated by three elements and then has the form $k[x,y,z]/(f)$. We have marked in boldface those canonical types, where $\La$ is minimal wild.
 \begin{center}
 \def\bu{\bullet}
 \begin{tabular}{|c|l|l|r|c|} \hline
 $\La$     & $\deg(x,y,z)$ & relation $f$    &$\deg{f}$    \\ \hline\hline
 \bf(2,3,7)& $(6,14,21)$   & $z^2+y^3+x^7$   & $42$       \\
 $(2,3,8)$ & $(6,8,15)$    & $z^2+x^5+xy^3$  & $30$       \\
 $(2,3,9)$ & $(6,8,9)$     & $y^3+xz^2+x^4$  & $24$        \\ \hline
 \bf(2,4,5)& $(4,10,15)$   & $z^2+y^3+x^5y$  & $30$        \\
 $(2,4,6)$ & $(4,6,11)$    & $z^2+x^4y+xy^3$ & $22$        \\
 $(2,4,7)$ & $(4,6,7)$     & $y^3+x^3y+xz^2$ & $18$        \\
 $(2,5,5)$ & $(4,5,10)$    & $z^2+y^2z+x^5$  & $20$  \\
 $(2,5,6)$ & $(4,5,6)$     & $xz^2+y^2z+x^4$ & $16$      \\ \hline
 \bf(3,3,4)& $(3,8,12)$    & $z^2+y^3+x^4z$  & $24$ \\
 $(3,3,5)$ & $(3,5,9)$     & $z^2+xy^3+x^3z$ & $18$  \\
 $(3,3,6)$ & $(3,5,6)$     & $y^3+x^3z+xz^2$ & $15$  \\
 $(3,4,4)$ & $(3,4,8)$     & $z^2-y^2z+x^4y$ & $16$  \\
 $(3,4,5)$ & $(3,4,5)$     & $x^3y+xz^2+y^2z$& $13$     \\
 $(4,4,4)$ & $(3,4,4)$     & $x^4-yz^2+y^2z$ & $12$  \\ \hline
 \end{tabular}
 \end{center}
For $k=\CC$ the equations are equivalent to \define{Arnold's exceptional unimodal singularities} in
the theory of singularities of differentiable maps~\cite{Arnold:1985}. In the theory of {\em
automorphic forms} the 14 graded algebras are known to be exactly those rings of entire
automorphic forms having \emph{three generators}~\cite{Wagreich:1980}.

\subsection{The triangulated category of singularities} \label{ssect:dsing}
For a variety $X$, Orlov investigated in \cite{Orlov:2004} the
triangulated category $\Dsing{X}$ of singularities of $X$ defined
as the quotient of the bounded derived category $\Der{\coh{X}}$ of
coherent sheaves modulo the full subcategory of perfect complexes.
If $X$ is affine with coordinate algebra $R$, then this category $\Dsing{R}$
is just the quotient $\Der{\mmod{R}}/\Der{\proj{R}}$, where
$\proj{R}$ is the category of finitely generated projective
$R$-modules. In \cite{Orlov:2005} Orlov further introduced a graded variant
$
\Dsinggrad\ZZ{R}=\Der{\modgrad\ZZ{R}}/\Der{\projgrad\ZZ{R}},
$
called the \emph{triangulated category of the graded singularities of
$R$} which will play a central role in this section. An `infinite' version of this category was studied in~\cite{Krause:2005}.

Under the name \emph{stabilized derived category of $R$} the
categories $\Dsing{R}$ were introduced before by Buchweitz in \cite{Buchweitz:1987}.
Extended to the graded case, Buchweitz gives an alternative description of $\Dsinggrad\ZZ{R}$ as
the \emph{stable category of graded maximal Cohen-Macaulay modules}
$\stabCMgrad{\ZZ}{R}$, in Orlov's treatment this category appears as the \emph{stable category of matrix factorizations}, used for instance in \cite{Kajiura:Saito:Takahashi:2006} and~\cite{Kajiura:Saito:Takahashi:2007}.  More precisely, Buchweitz showed that the category
$\CMgrad{\ZZ}{R}$ of maximal graded Cohen-Macaulay $R$-modules is a
Frobenius-category, hence inducing,
--- in Keller's terminology \cite{Keller:2006} --- on the attached stable
category $\stabCMgrad{\ZZ}{R}$ of graded maximal Cohen-Macaulay
modules modulo projectives,  the structure of an \emph{algebraic}
triangulated category, see \cite{Happel:1988} for the definition of a Frobenius category and its attached stable triangulated category.

\subsection{Orlov's theorem} \label{ssect:orlov}
Orlov's theorem states an interesting trichotomy relating the derived category $\Dd=\Der\Cc$ of coherent sheaves on a \define{generalized projective space}, given by a positively $\ZZ$-graded coordinate algebra $A$, to the triangulated category $\Dsinggrad\ZZ{A}$ of singularities of $A$. Under suitable restrictions on $A$, the category of coherent sheaves $\Cc$ in question is defined as $\modgrad\ZZ{A}/\modgradnull\ZZ{A}$.

\begin{theorem}[Orlov~\cite{Orlov:2005}]\label{thm:orlov}
Let $A$ be a positively $\ZZ$-graded noetherian algebra which is graded Gorenstein with Gorenstein index $a$. Assume that $A$ is connected, that is, $A_0=k$. Then the triangulated categories $\Dsinggrad\ZZ{A}$ and $\Dd=\Der{\modgrad\ZZ{A}}/\modgradnull\ZZ{A}$ are related as follows:

$(i)$ If $a>0$, then $\Dsinggrad\ZZ{A}$ can be identified with a subcategory of $\Dd$ which is the (right) perpendicular category to an exceptional sequence of length $a$.

$(ii)$ If $a=0$, there exists a triangle-equivalence $\Dsinggrad\ZZ\cong \Dd$.

$(iii)$ If $a<0$, then $\Dd$ can be identified with a full subcategory of $\Dsinggrad\ZZ{A}$ which is the (right) perpendicular category to an exceptional sequence of length $|a|$.~\hfill\qed
\end{theorem}
We are not giving the formal definition of the Gorenstein index $a$ in general. Instead we remark that, for a weighted projective line $\XX$, the Gorenstein index of $S=\Spla$ is zero in the tubular case, for $R=\Rpla$ it is $+1$, if $\eulerchar\XX$ is positive, and it is $-1$ if $\eulerchar\XX$ is negative. (The claim follows from Theorem~\ref{thm:gorenstein}.)

The categories of coherent sheaves on a weighted projective line offer an interesting application of Orlov's theorem. With some extra work, concerning the exact shape of the mentioned exceptional sequence and its perpendicular category, one gets the following result which is due to several authors \cite{Kajiura:Saito:Takahashi:2006},\cite{Kajiura:Saito:Takahashi:2007}, \cite{Ueda:2006} and \cite{Lenzing:Pena:2006}. Note that the authors of \cite{Kajiura:Saito:Takahashi:2006,Kajiura:Saito:Takahashi:2007} work in the stable category of matrix factorizations, so in a copy of $\stabCMgrad\ZZ{A}$, instead of $\Dsinggrad\ZZ{A}$. In the tubular case, one uses a natural generalization of Orlov's theorem to allow $\Lp$-gradings, see~\cite{Ueda:2006}. The next (sub)section will contain a sketch of the proof of Orlov's theorem in the case of weighted projective lines.

\subsection{An application of Orlov's theorem} \label{ssect:appl_orlov}
The next result illustrates the impact of Orlov's theorem. It is not a simple corollary, however, since serious extra work has to be done to get the explicit shape of the result.
\begin{theorem} \label{thm:trichotomy}
  Let $\XX$ be a weighted projective line. Then the following holds:

$(i)$ If $\eulerchar\XX>0$, let $R=\Rpla$ be the attached $\ZZ$-graded singularity and $\vDe=[p_1,p_2,p_3]$ be the Dynkin star, obtained from the weight sequence of $\XX$. Then there is a triangle-equivalence $\Dsinggrad\ZZ{R}\cong\Der{\mmod{k[\vDe]}}$.

$(ii)$ If $\eulerchar\XX=0$, let $S=\Spla$ be the attached $\Lp$-graded singularity. Then there is a triangle-equivalence $\Dsinggrad\Lp{S}\cong\Der{\coh\XX}$.

$(iii)$ If $\eulerchar\XX<0$, let $R=\Rpla$ be the attached $\ZZ$-graded singularity and $\Lahat=\Lahat\spitz{\pla}$ be the extended canonical algebra attached to $\XX$. Then there is a triangle-equivalence $\Dsinggrad\ZZ{R}\cong\Der{\mmod{\Lahat}}$.
\end{theorem}
Concerning $(i)$ we refer to \cite{Kajiura:Saito:Takahashi:2006}, also a proof in the discussed setting would not be difficult.
As observed in~\cite{Ueda:2006}, assertion $(ii)$ follows directly from an $\Lp$-graded version of Orlov's theorem. \emph{From now on we are dealing with case $(iii)$.} We need some preparation for the proof.

Let $R=\Rpla$, be the positively
$\ZZ$-graded Gorenstein singularity attached to $\XX$. As mentioned before, $R$ has Gorenstein index $-1$. We fix some notation:
Let $\Mm=\Der{\modgrad{\ZZ}{R}}$ and $\Mmplus=\Der{\modgrad\ZZplus{R}}$. Let $\Ppplus$ be the triangulated subcategory of $\Mmplus$
generated by all $R(-n)$, $n\geq0$ and $\Tt$ its left perpendicular
category $\lperp{\Ppplus}$ formed in $\Mmplus$. Denote further by
$\Ssplus$ the triangulated subcategory of $\Mmplus$ generated by all
$k(-n)$, $n\geq0$ and $\Dd$ its right perpendicular category
$\rperp{\Ssplus}$ formed in $\Mmplus$. Finally let
$\Dd(-1)=\rperp{\Ssplus(-1)}$. From the proof of Orlov's theorem \cite[2.5]{Orlov:2005} we deduce the
following proposition.

\begin{proposition} \label{propo:special}
 Assume that $\eulerchar{\XX}<0$ and let $R$ be the
positively $\ZZ$-graded singularity attached to $\XX$. Then the
following holds:

$\mathrm{(a)}$ The natural functor $\Tt \incl \Mm\stackrel{q}{\ra}
\Dsing{R}$, where $q$ is the quotient functor, is an equivalence of
triangulated categories.

$\mathrm{(b)}$ The $R$-module $k$ is an exceptional object in $\Tt$
with $\lperp{k}=\Dd(-1)$. Moreover, the category $\Der{\coh\XX}$
is naturally equivalent to $\Dd(-1)$ under the functor $Y\mapsto
(\bR \Gaplus(Y))(-1)$.
\end{proposition}
\begin{proof}
We sketch the
argument: Using that $R$ is Gorenstein of Gorenstein index -1, and
invoking Gorenstein duality $\bR\Homgrad{\bullet}{R}{-}{R}$ of $\Mm$ one
sees that $\rperp{T}\subset \rperp{\Dd(-1)}$ and hence $\Dd(-1)$ is
a full subcategory of $\Tt$. Further we see that
$\lperp{k}=\Dd(-1)$. It is well-known that
$
\Gaplus:\coh\XX \ra \modgrad\ZZplus{R}$,
$Y\mapsto
\bigoplus_{n=0}^\infty \Hom{}{\Oo}{Y(n)}$,
is a full embedding having sheafification, that is, the quotient
functor $q_+:\modgrad\ZZplus{R}\ra \Dsinggrad\ZZ{R}$ has an exact left
adjoint and such that composition $q\,\Gaplus$ is the identity
functor on $\Dsinggrad\ZZ{R}$, compare~\cite[1.8,5.1]{Geigle:Lenzing:1987} and \cite[5.7]{Lenzing:1994}. It
follows that $\bR\Gaplus:\Der{\coh\XX}\ra \Mmplus$ is a full
embedding having $q_+:\Mmplus\ra\Mmplus/\Ssplus$ as a left adjoint,
and $q_+\,\bR\Gaplus=1$.

Since $R$ is positively graded with $R_0=k$, it follows that $k$ is
exceptional in $\Mm$ and hence in $\Mmplus$. Invoking the minimal
graded injective resolution $0\ra R \ra E^0\ra E^1 \ra E^2\ra 0$,
where $E^0$ and $E^1$ are socle-free and $E^2$ is the graded
injective hull of $k(-1)$, it follows that $k$ belongs to $\Tt$ and
then also to $\Dd(-1)$. It is straightforward to check that
$\rperp{\Dd(-1)}$ equals the triangulated subcategory $\spitz{k}$ generated by $k$, and hence $\lperp{k}=\Dd(-1)$ in $\Tt$.
\end{proof}
We now return to the \textbf{Proof of $(iii)$.} The idea is: $(a)$ to take a tilting object $T$ in $\coh\XX$, whose endomorphism ring is derived equivalent to a canonical algebra $\La$ attached to $\XX$ and, identifying $\coh\XX$ and $\Dd(-1)$, $(b)$ to show that $\hat{T}=T(-1)\oplus k$ yields a tilting object in $\Tt$, and $(c)$ to show that the endomorphism ring of $\hat{T}$ is derived equivalent to the extended canonical algebra, attached to $\XX$. To obtain $(b)$, it suffices by Proposition~\ref{prop:one-point} to show that the left
adjoint $\ell:\Tt\ra \lperp{k}$ to the inclusion $j:\lperp{k}\incl
\Tt$ maps $k$ to a line bundle in $\Dd(-1)=\Der{\coh\XX}$ up to
translation in $\Dd(-1)$. We put $A=(\bR\Gaplus(\Oo(-\vom)))(-1)$
and construct a morphism $\ga:k\ra A[1]$ such that
$\Hom{}{\ga}{Y}:\Hom{}{A[1]}{Y}\ra\Hom{}{k}{Y}$ is an isomorphism for each
$Y\in\lperp{k}$ such that $\ell(k)=A[1]$.

The claim is proved in two steps. Put $R_+=\bigoplus_{n\geq 1}R_n$,
then the exact sequence $0\ra R_+\ra R \ra k$ yields an exact
triangle $R\ra k \stackrel{\al}{\ra} R_+[1]$ in $\Mm$, where
$\Hom{}{\al}{Y}$ is an isomorphism for each $Y\in\lperp{k}$. Note for
this that $R$ belongs to $\lperp{\Dd(-1)}$.

For the next step it is useful to identify the derived category
$\Mmplus$ with the full subcategory of $\Der{\Modgrad\ZZplus{R}}$
consisting of all complexes with cohomology in $\modgrad\ZZplus{R}$. Here, $\Modgrad\ZZplus{R}$ denotes the category of \emph{all}
graded $R$-modules. Let $0\ra R(-1)\ra E^0 \ra E^1 \ra E^2 \ra 0$ be
the minimal graded injective resolution of $R(1)$ such that $E^2$
equals the graded injective envelope of $k$. (This uses that $R$ has
Gorenstein index $-1$.) Sheafification yields the minimal injective
resolution $0\ra \Oo(\vom)\ra \tilde{E}^0\ra \tilde{E}^1 \ra0$ of
$\Oo(\vom)$. Accordingly $\bR\Gaplus(\Oo(\vom))$ is given by the
complex
$
A: \quad \cdots\ra 0 \ra E_+^0(-1) \ra E_+^2(-1) \ra 0 \cdots,
\textrm{ where } E_+^i=\bigoplus_{n\geq0}E_n^i,
$
whose cohomology is concentrated in degrees zero and one and given
by
$
\mathrm{H}^0(A)=R_+$ and $ \mathrm{H}^1(A)=k(-1).
$
It follows the existence of an exact triangle
$
k(-1)[-2]\ra R_+ \stackrel{\be}{\ra} A \ra k(-1)[-1],
$
in $\Mmplus$ where, by construction, $A$ belongs to $\Dd(-1)$. For
$Y$ from $\Dd(-1)$ we have, in particular, that $Y$ belongs to
$\lperp{k(-1)}$ implying that $\Hom{}{\be}{Y}$ is an isomorphism. To
summarize: The morphism $\ga=[k\stackrel{\al}{\ra} R_+[1]
\stackrel{\be[1]}{\ra} A[1]]$ yields isomorphisms $\Hom{\ga}{Y}$ for
each $Y\in \Dd(-1)$. Hence $\ell(k)[-1]=A=\bR\Gaplus(\Oo(\vom))(-1)$
is a line bundle, as claimed. The claim now follows from Proposition~\ref{prop:special}.\qed

\section{The stable category of vector bundles}\label{sect:stable}
Let $\XX$ be a weighted projective line. Relative to the natural choice of a class of $\Ll$ of line bundles on $\XX$ we give the category $\vect\XX$ of vector bundles on $\XX$ the structure of a {Frobenius category} whose indecomposable projectives=injectives are just the members of $\Ll$. We recall that a $k$-linear category $\Ff$ is called a \define{Frobenius category} if it is equipped with an \emph{exact structure} inherited from an abelian category through a full embedding, if further $\Ff$ has enough projectives and injectives with respect to this structure, and finally projectives and injectives coincide in $\Ff$. By~\cite{Happel:1988} the corresponding stable category $\underline{\Ff}$, that is,  the factor category of $\Ff$ by the ideal $\Ii$ of all morphism factoring through a projective (or injective), is a triangulated category (which is algebraic in the sense of~\cite{Keller:2006}).

\subsection{Fundamental properties} \label{ssect:fundamental_prop}
For this section we use unpublished material from joint work with D.~Kussin and H.~Meltzer~\cite{Kussin:Lenzing:Meltzer:Pena:2008}. In this account we concentrate on the stable category of vector bundles attached to a weighted projective line, although a more general approach is possible.

Let $\XX$ be a weighted projective line. We define the \define{distinguished class of line bundles} $\Ll$ on $\XX$ as follows: If $\eulerchar\XX\neq 0$, let $\Ll=\tau$-orbit $\{\Oo(n\vom)| n\in\ZZ \}$ is the $\tau$-orbit of the structure sheaf. In the tubular case $\eulerchar\XX=0$ we take the system $\Ll=\{\Oo(\vx)|\vx\in\Lp\}$ of all line bundles. Note that in each case $\Ll$ is closed under the Auslander-Reiten translation. It is moreover convenient, to assume that the class of distinguished line bundles is closed under isomorphism. (Further choices for a distinguished class are possible and discussed in \cite{Kussin:Lenzing:Meltzer:Pena:2008}). By definition, the \define{stable category of vector bundles}\footnote{This should not be confused with the category of stable vector bundles!} $\stabvect{\XX}$ is the triangulated category obtained from $\vect\XX$ as the factor category $\vect\XX/[\Ll]$ of $\vect\XX$ modulo the two-sided ideal of all morphisms factoring through a finite direct sum of members from $\Ll$. We are going to show that this stable category has a natural triangulated structure.

\begin{definition}
A sequence $\eta:0\ra E'\ra E \ra E'' \ra 0$ in $\vect\XX$ is called a \define{distinguished exact sequence} in $\vect\XX$ if the functor $\Hom{}{L}{-}$ is exact on $\eta$ for each distinguished line bundle $L$.
\end{definition}

It follows from \cite{Buchweitz:1987} that the class of sequences in $\CMgrad\ZZ{R}$, resp.\ $\CMgrad\Lp{S}$ which are exact in the category of all graded modules over $R=\Rpla$ (resp.\ $S=\Spla$) is the class of distinguished exact sequence for the exact structure of a Frobenius category on the category of graded maximal Cohen-Macaulay modules in question. By Theorem~\ref{thm:gorenstein} we have natural identifications $\CMgrad\ZZ{R}=\vect\XX$ for $\eulerchar\XX\neq0$ (resp.\ $\CMgrad\Lp{S}=\vect\XX$) for $\eulerchar\XX=0$ identifying indecomposable projective (graded) CM-modules with distinguished line bundles. Moreover, these identifications induce an exact structure on $\vect\XX$ which agrees with the above exact structure. We thus obtain:
\begin{theorem} \label{thm:frobenius}
The class of distinguished exact sequences defines on $\vect\XX$ an exact structure turning $\vect\XX$ into a Frobenius category whose indecomposable projectives (injectives) form the class $\Ll$ of distinguished line bundles.
Each distinguished exact sequence in $\vect\XX$ is an exact sequence in the abelian category $\coh\XX$, but not conversely. Moreover there are natural identifications $\stabCMgrad\ZZ{R}=\stabvect\XX$ for $\eulerchar\XX\neq0$ and $\stabCMgrad\Lp{S}=\stabvect\XX$ for $\eulerchar\XX=0$.~\hfill\qed
\end{theorem}

\begin{theorem} \label{thm:stable_cat}
Let $\XX$ be a weighted projective line of weight type $\pP=(p_1,\ldots,p_t)$. Then the category $\stabvect\XX$ is a triangulated category with Serre duality which has a tilting object. Moreover a tilting object $T$ can be chosen in such a way such that:

$(i)$ $\End{}{T}$ is the path algebra of the Dynkin quiver $[p_1,\ldots,p_t]$ if $\eulerchar\XX>0$.

$(ii)$ $\End{}{T}$ is the canonical algebra of (tubular) weight type $\pP$ if $\eulerchar\XX=0$.

$(iii)$ $\End{}{T}$ is an extended canonical algebra of type $\spitz{p_1,\ldots,p_t}$ if $\eulerchar\XX<0$.
\end{theorem}
\begin{proof}
Use Theorems~\ref{thm:frobenius} and \ref{thm:trichotomy}.
\end{proof}

Note that the algebras $\End{}{T}$ above all have finite global dimension. From the triangle equivalence $\stabvect\XX\iso\Der{\mmod{\End{}{T}}}$ we therefore conclude.

\begin{coro}
The category $\stabvect\XX$ is homologically finite and has Serre duality. Moreover the Grothendieck group of $\stabvect\XX$ (viewed as a  triangulated category) is finitely generated free and the Euler form on $\stabvect\XX$, given on classes of objects by
 $$
 \eulerform{[X]}{[Y]}=\sum_{n\in\ZZ}(-1)^n\ddim{k}{\Hom{}{X}{Y[n]}},
 $$
is non-degenerate.~\hfill\qed
\end{coro}

For the proof of the next proposition we refer to \cite{Kussin:Lenzing:Meltzer:Pena:2008}.

\begin{proposition}
$(i)$ Let $X$ and $Y$ be indecomposable  vector bundles not from $\Ll$. Then a morphism in $u:X\ra Y$ is irreducible in $\vect\XX$ if and only if the induced morphism $\ul{u}:\ul{X}\ra \ul{Y}$ is irreducible in $\stabvect\XX$.

$(ii)$ Each Auslander-Reiten sequence $0\ra E' \ra E \ra E'' \ra 0$ in $\vect\XX$ whose end terms are not in $\Ll$ is a distinguished exact sequence and yields a distinguished triangle in $\stabvect\XX$.

$(iii)$ The Auslander-Reiten translation on $\stabvect\XX$, as a functor, is induced by the Auslander-Reiten translation on $\vect\XX$.
\end{proposition}

\begin{coro}
Each Auslander-Reiten component $\Cc$ of $\vect\XX$ yields an Auslander-Reiten component $\ul\Cc$ of $\stabvect\XX$, and each Auslander-Reiten component of $\stabvect\XX$ is obtained this way. In more detail:

$(i)$ If $\eulerchar\XX>0$ with attached Dynkin diagram $[p_1,\ldots,p_t]$, then the indecomposables of $\stabvect\XX$ form a unique Auslander-Reiten component having type $\ZZ\De$.

$(ii)$ If $\eulerchar\XX=0$, then the components of $\stabvect\XX$ form a rational family (indexed by $\QQ \union \{\infty\}$) of one-parameter family of tubes $\Tt_\la$,  $\la\in \Proj{1}{k}$.

$(iii)$ If $\eulerchar\XX<0$, then each component of the Auslander-Reiten quiver of $\stabvect\XX$ has the form $\ZZ\AA_\infty$.~\hfill\qed
\end{coro}
On a theoretical level, the results of this section look very complete. However, several questions remain: \emph{How can we explicitly construct tilting objects in $\stabvect\XX$?} \emph{How do the concepts ``exceptional bundle in $\vect\XX$'' and ``exceptional object in $\stabvect\XX$'' relate?}
\subsection{Positive Euler characteristic} \label{ssect:euler_positive}
For $\eulerchar\XX>0$ it is easy to specify a tilting object. Let $p$ be the least common multiple of the weight sequence. The degree is the additive function on $\coh\XX$ which is zero on the structure sheaf, which is one on a simple sheaf concentrated in an ordinary point of $\XX$, and finally is  $p/p_i$ on each simple sheaf concentrated in an exceptional point of weight $p_i$. The slope $\mu X$ is the quotient of degree and rank.
\begin{proposition} Assume $\eulerchar\XX>0$. Let $\De=[p_1,p_2,p_3]$ be the Dynkin diagram given by the weights of $\XX$. Then the following holds:

\emph{(i)} The isomorphism classes of indecomposable objects in $\coh\XX$ with a slope in the range $0\leq q <p\eulerchar\XX$ form a finite system $\Ee$. The direct sum $T$ of all objects in $\Ee$  is a tilting object for the \emph{abelian category} $\coh\XX$, and $H=\End{}{T}$ is the path algebra of  extended Dynkin type $\tilde\De$.

\emph{(ii)} The direct sum of all objects in  $\Ee'=\Ee\setminus\{\Oo\}$ is a tilting object for the \emph{triangulated category} $\Tt=\stabvect\XX$, and $\End\Tt{T'}$ is a path algebra of Dynkin type $\De$.
\end{proposition}
\begin{proof}
The first claim is proved in~\cite{Lenzing:Reiten:2006}. It is then straightforward to prove the second claim. (Note however that (ii) is not just a trivial consequence of (i).)
\end{proof}
Note that this reestablishes the central result of \cite{Kajiura:Saito:Takahashi:2006} bypassing Orlov's theorem.

\subsection{Euler characteristic zero} \label{ssect:euler_zero}
As in \cite{Lenzing:Meltzer:1993} let $\Cc^{(q)}$ be the additive closure of all indecomposables from $\coh\XX$ with slope $q$, where $q\in\QQ\union\{\ \infty\}$. There is a self-equivalence $\varphi$ of $\Der{\coh\XX}$ sending $\Cc^{q}$ to $\Cc^{q/(1+q)}$ for each $q$. The following result, taken from~\cite{Kussin:Lenzing:Meltzer:Pena:2008}, has a straightforward proof.

\begin{proposition}
Assume that $\eulerchar\XX=0$.
Then the direct sum of all $\varphi(\Oo(\vx))$, with $\vx\in\Lp$ such that $0\leq \vx \leq \vc$, is a tilting object in $\stabvect\XX$ whose endomorphism ring is the canonical algebra $\La$ attached to $\XX$.~\hfill\qed
\end{proposition}
\begin{coro}
The triangulated categories $\stabvect\XX$ and $\Der{\coh\XX}$ are triangle-equivalent.~\hfill\qed
\end{coro}
This gives a direct proof of Ueda's result~\cite{Ueda:2006} not relying on Orlov's theorem.

\subsection{Negative Euler characteristic: weight type $(2,3,7)$} \label{ssect:euler_negative}
We next deal with the shape of $\stabvect\XX$ for Euler characteristic $\eulerchar\XX<0$. By the preceding we know the combinatorial structure of $\stabvect\XX$, in particular that each Auslander-Reiten component has shape $\ZZ\AA_\infty$. Little information, we have sofar concerning the categorical structure of these components and of the morphisms between different components. In particular, it is difficult to give the precise location of a tilting object whose endomorphism ring is extended canonical.

Fairly good information is available, however, on those components of $\stabvect\XX$ arising from components containing a line bundle (and then a whole Auslander-Reiten orbit of line bundles). We call them \define{line bundle components}\footnote{This, despite the fact that in the stable category the line bundles may no longer be visible.}. The number of line bundle components
 $
 p_1\cdots p_t\left((t-2)-\sum_{i=1}^t{1/p_i}\right)
 $
is always finite, and only depends on the weight sequence~\cite[prop.~8.2]{Lenzing:Pena:1997}. This number is minimal for weight type $(2,3,7)$ where the line bundles form a single $\tau$-orbit, and $(2,3,7)$ is the only weight type where this happens.

We are going to exhibit a tilting object for $\stabvect\XX$ in this line bundle component $\Dd$ such that, moreover, its indecomposable summands are from a single Auslander-Reiten orbit.

\begin{theorem} \label{thm:tilt:orbit}
We assume that $\XX$ has weight type $(2,3,7)$. Let $M$ be the set of all integers of the form $22a+7b$ with $a\in\set{0,\ldots,5}$ and $b\in\set{0,1}$. Let $E$ be the extension term of the almost-split sequence $0\ra \Oo\ra E \ra \tau^{-1}\Oo \ra 0$. Then $T=\Dir_{n\in M}\tau^n E$ is a tilting object in $\stabvect\XX$, whose endomorphism algebra is given by the linear quiver of $12$ vertices with all zero relations $x^3=0$.
\end{theorem}

The proof affords several steps. We start with the key result.

\begin{proposition} \label{prop:key}
We have $\ulHom{}{E}{\tau^n E}\neq0$ for an integer $n$ if and only if $n$ equals $0$, $7$ or $22$.
\end{proposition}
\begin{proof}
We abbreviate $\tau^nF$ by $F(n)$. Then
$$
0\ra E(-21)\up{} \Oo(-20)\ddir\Oo(-15)\ddir\Oo(-7)\ddir\Oo\up{} E \ra0
$$
is the minimal projective resolution of $E$ in the Frobenius category $\vect\XX$. The claim follows by straightforward computation.
\end{proof}
\begin{coro}
$E$ is quasi-simple in the stable line bundle component $\ul\Dd$. All objects having quasi-length $\leq 6$ in $\ul\Dd$ are exceptional. Further, the objects of quasi-length $7$ have trivial endomorphism ring $k$, but are not exceptional.
\end{coro}
\begin{proof}
By the shape of the Auslander-Reiten component, up to repeated $\tau$-shift, an object $F$ of quasi-length $n$ has a \define{triangle filtration} $0=F_0\ra F_1 \ra \cdots \ra F_n=F$ such that for each $j=1,\ldots,t$ we have a triangle $F_{j-1}\ra F_{j}\ra \Oo(n-j)\ra F_{j-1}[1]$. The assertion now immediately follows from the proposition.
\end{proof}
Calculating minimal projective resolutions in $\vect\XX$ we obtain:
\begin{lemma}\label{lemma:tau_hoch_21}
On the stable line bundle component $\ul\Dd$ the translation $[1]$ acts as $\tau^{21}$.~\hfill\qed
\end{lemma}

\begin{proposition} \label{prop:extensionless}
With the notations above, we have $\ulHom{}{E(\ell)}{E(\ell'[n])}=0$ for all $\ell,\ell'$ from $M$ and all integers $n\neq0$.
\end{proposition}
\begin{proof}
We assume that $\ulHom{}{E(\ell)}{E(\ell'[n])}\neq0$ and
write $\ell=22a+7b$ and $\ell'=22a'+7b'$ with $a,a'\in\set{0,\ldots,5}$ and $b,b'\in\set{0,1}$. Then $$
\alpha:=22(a'-a)+7(b'-b)+21n=\ell'-\ell+21n\in\set{0,7,22}.
$$
There are two cases to consider: \ul{Case 1} $\alpha\in\set{0,7}$. Then $a'-a=0 \mmod 7$, and $a'=a$ follows. Moreover, $b'-b+3n$ belongs to $\set{0,1}$ which implies $n=0$, and proves the claim in this case. \ul{case 2} $\alpha=22$. Reduction modulo $7$ yields $a'-a=1\mmod 7$, hence $a'=a+1$ and $(b'-b)+3n=0$ which in turn implies $n=0$ also in this case.
\end{proof}
It affords more work to show that the system $\set{M(\ell)|\ell\in M}$ generates $\stabvect\XX$ as a triangulated category.

We start with a result that is independent of the weight type.

\begin{proposition}
Let $X$ be an indecomposable vector bundle. Let $L$ be a line bundle of maximal degree such that $\Hom{}{L}{X}\neq0$. Each non-zero $f:L\ra X$ extends to the extension term $E$ of the almost-split sequence $0\ra L\ra E \ra \tau^{-1}L \ra 0$ and yields a morphism $\bar{f}:E\ra X$ that does not factor through a direct sum of line bundles.~\hfill\qed
\end{proposition}

\begin{proposition} \label{prop:generate}
Let $T=\Dir_{\ell\in M}E(\ell)$ and $X$ an indecomposable vector bundle not in $\Ll$. Then $\ulHom{}{T}{X[n]}\neq0$ for some integer $n$.~\hfill\qed
\end{proposition}
\begin{proof}
By the preceding proposition, we know that $\ulHom{}{E(\ell)}{X}\neq0$ for some integer $\ell\in\ZZ$. By Lemma~\ref{lemma:tau_hoch_21} we can reduce $\ell$ modulo 21. Note that reduction modulo 21 for $\Mm$ yields the set $S=\{0,1,\ldots,5,7,\ldots,12\}$. For this we use the exact triangle $E(0)\up{\bar{x_1}}E(7)\up{\bar{x_1}}E(14)\up{x_1}E(21)$, where $\bar{x_1}$ denotes the multiplication by $x_1$. This allows to enlarge $S$ by $\{14,15,16,17,18,19\}$. To obtain the missing elements 6,13,20 we establish a further exact triangle coming from a distinguished exact sequence $0\ra E(5)\ra \Oo(4)\oplus E(12)\oplus E(20)\ra F\ra 0$, where $F$ has quasi-length $3$ and quasi-socle $\Oo(20)$ in the line bundle component of $\coh\XX$, compare the following piece of the Auslander-Reiten quiver in $\coh\XX$:
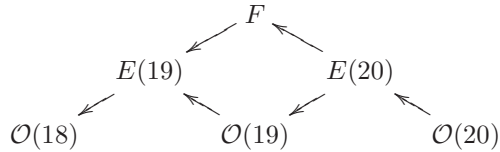
$$
\xymatrix@R8pt@C8pt{
      &         &F\ar[ld]&        &         \\
      &E(19)\ar[ld]&        &E(20)\ar[lu]\ar[ld]&        \\
\Oo(18)&        &\Oo(19)\ar[lu] &         & \Oo(20)\ar[lu]\\
}
$$
We have shown that $\ulHom{}{E}{X}\neq0$.
\end{proof}

\subsubsection*{Proof of theorem~\ref{thm:tilt:orbit}} It follows from propositions~\ref{prop:extensionless} and \ref{prop:generate} that $T$, the direct sum of all bundles $E(\ell)$ with $\ell\in M$, is a tilting object in the triangulated category $\stabvect\XX$. By means of proposition~\ref{prop:key} it is further easily checked that the endomorphism algebra of $T$ is isomorphic to the algebra $A_{12}$ from Section~\ref{sect:accessible} given by the linear quiver of $12$ vertices, satisfying all zero-relations $x^3=0$.~\hfill\qed

\section{Discussion and comments} \label{sect:comments}

\subsection{The Coxeter formalism}
For a finite dimensional algebra $A$ of finite global dimension the \define{Coxeter formalism} given by (a) the Euler bilinear form on $\Knull{A}$, (b) the Coxeter transformation induced from the Auslander-Reiten translation $\tau$ of $\Der{\mmod{A}}$, and (c) the Coxeter polynomial $\coxpol{A}$,  yields a particularly convincing set of data reflecting deep homological properties of $A$.
There are other contexts in mathematics where a Coxeter formalism exists. Such connections often indicate a relationship to the representation theory of finite dimensional algebras. Such links therefore should be followed in order to determine whether the formal relationship is based on a conceptual connection.

(i) We have already mentioned \emph{graph theory} in Section~\ref{ssect:graph}.

(ii) There is a now classical relationship to \emph{Lie theory}, where the Coxeter formalism originates. The first link is the occurrence of Dynkin diagrams: By Gabriel's theorem~\cite{Gabriel:1972} the isomorphism classes of indecomposable representations of the path algebra $k[\vDe]$ of a Dynkin quiver are in bijective correspondence to their classes in the Grothendieck group $\Knull{k[\vDe]}$, where they form the root system given by the roots of the corresponding \emph{Tits} or \emph{Euler quadratic form}. The correspondence in the converse direction is due to Ringel, who attaches a (quantum) Lie algebra to $k[\vDe]$. The basic construction is now known under the name \emph{Ringel-Hall algebra}~\cite{Ringel:1994a}.

(iii) Another instance is \emph{knot theory}. There `the' \define{Seifert matrix} $V$ takes the role of the Cartan matrix, and `the' \define{Alexander polynomial} $\det(V^t-xV)$ takes the role of the Coxeter polynomial, up to the factor $\det(V)$. We refer to \cite{Burde:Zieschang:2003} for the relevant definitions.

(iv) This survey is devoted to the link to \emph{singularity theory}. We assume $k=\CC$ and briefly review some concepts; the reader will find all relevant definitions in~\cite{Ebeling:2007}. For the two-dimensional singularities, appearing in this paper, the \emph{Milnor lattice} is the second integral homology $\tief{2}{X}$ of the \emph{Milnor fibre}, equipped with the symmetric bilinear \emph{intersection form}. With respect to a \emph{strongly distinguished basis} the \emph{variation matrix} $V$, the \emph{intersection matrix} $S$ and the \emph{classical monodromy} $H$ are related by
$$
S=-(V+V^t),\quad \textrm{and } H=-V^{-1}V^t,
$$
showing the existence of a Coxeter formalism. Relative to a \emph{strongly distinguished basis}, the intersection matrix $S$ further yields a bigraph (a graph with two kind of edges, say solid and dotted), called a \emph{Coxeter-Dynkin diagram} of the singularity. We note that for a derived accessible algebra of simpler type (Dynkin, extended Dynkin, canonical, extended canonical) the underlying quiver yields the solid edges of an associated Coxeter-Dynkin diagram, where the dotted lines represent the relations. We refer to \cite{Brieskorn:1987} and \cite{Ebeling:96} for a detailed investigation of Milnor lattices.

In conclusion, our approach yields a categorification of the Milnor lattice for the singularity attached to a weighted projective line by a triangulated category with a tilting object, that is, by a bounded derived category of a certain finite dimensional algebra $A$. Note, that for the singularities, discussed in this survey, the algebra $A$ is derived accessible and, moreover, strongly related to hereditary representation theory.

For Euler characteristic $\eulerchar\XX<0$ an interesting twist occurs: For the weight types $(p,q,r)$ from Arnold's list of exceptional unimodal singularities the Coxeter-Dynkin algebra we attach to $\XX$, see Section~\ref{ssect:ext_canonical}, yields `the' Coxeter-Dynkin diagram of the Milnor lattice corresponding to a singularity $R(p',q',r')$, where the mapping $(p,q,r)\mapsto (p',q',r')$ is known as \emph{Arnold's strange duality}, an instance of mirror symmetry, see~\cite{Kobayashi:1995}

\subsection{The role of derived accessible algebras}
The following picture yields a rough description of the panorama of (derived accessible) algebras. This picture is well understood as long as we keep close to the center, given by the hereditary algebras of Dynkin or extended Dynkin type. It is not so clear what happens if we continue to extend these algebras by exceptional modules, when starting from wild hereditary or wild canonical algebras or even going beyond that. The knowledge is going to get very poor, when we start to extend supercanonical or extended canonical algebras. This is due to the insufficient knowledge of exceptional modules and the structure of their one-point extensions. Further work has to be done in this context. Here, spectral analysis has a pilot function, because once the dimension vector (class in the Grothendieck group) of an exceptional module is known, then at least the Coxeter polynomial of the one-point extension is available.
 \begin{center}
\setlength{\unitlength}{2700sp}%
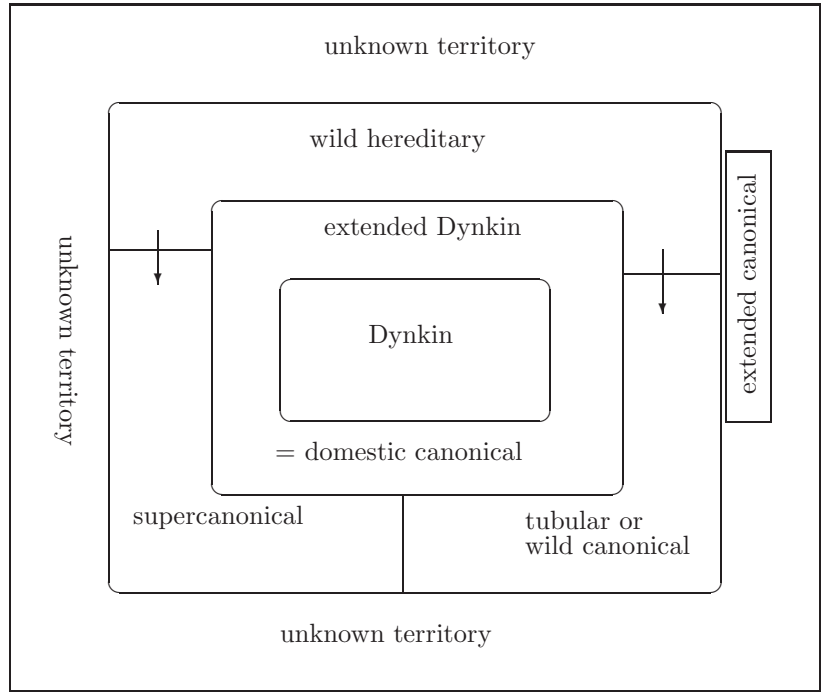
\begin{picture}(7449,6324)(2239,-6823)
\thinlines
{\put(4831,-4231){\oval(210,210)[bl]}
\put(4831,-3136){\oval(210,210)[tl]}
\put(7096,-4231){\oval(210,210)[br]}
\put(7096,-3136){\oval(210,210)[tr]}
\put(4831,-4336){\line( 1, 0){2265}}
\put(4831,-3031){\line( 1, 0){2265}}
\put(4726,-4231){\line( 0, 1){1095}}
\put(7201,-4231){\line( 0, 1){1095}}
}%
{ \put(4201,-4906){\oval(210,210)[bl]}
\put(4201,-2416){\oval(210,210)[tl]}
\put(7771,-4906){\oval(210,210)[br]}
\put(7771,-2416){\oval(210,210)[tr]}
\put(4201,-5011){\line( 1, 0){3570}}
\put(4201,-2311){\line( 1, 0){3570}}
\put(4096,-4906){\line( 0, 1){2490}}
\put(7876,-4906){\line( 0, 1){2490}}
}%
{ \put(3256,-5806){\oval(210,210)[bl]}
\put(3256,-1516){\oval(210,210)[tl]}
\put(8671,-5806){\oval(210,210)[br]}
\put(8671,-1516){\oval(210,210)[tr]}
\put(3256,-5911){\line( 1, 0){5415}}
\put(3256,-1411){\line( 1, 0){5415}}
\put(3151,-5806){\line( 0, 1){4290}}
\put(8776,-5806){\line( 0, 1){4290}}
}%
{ \put(7876,-2986){\line( 1, 0){900}}
}%
{ \put(5851,-5011){\line( 0,-1){900}}
}%
{ \put(3151,-2761){\line( 1, 0){945}}
}%
{ \put(8236,-2761){\vector( 0,-1){585}}
}%
{ \put(3601,-2581){\vector( 0,-1){495}}
}%
{ \put(2251,-6811){\framebox(7425,6300){}}
}%
{ \put(8821,-4336){\framebox(405,2475){}}
}%
\put(5536,-3661){\makebox(0,0)[lb]{Dynkin}%
}
\put(5131,-2671){\makebox(0,0)[lb]{extended Dynkin}%
}
\put(4681,-4696){\makebox(0,0)[lb]{= domestic canonical}%
}
\put(4996,-1861){\makebox(0,0)[lb]{wild hereditary}%
}
\put(3376,-5326){\makebox(0,0)[lb]{supercanonical}%
}
\put(6976,-5326){\makebox(0,0)[lb]{tubular or}%
}
\put(6976,-5551){\makebox(0,0)[lb]{wild canonical}%
}
\put(9091,-4080){\rotatebox{90.0}{\makebox(0,0)[lb]{extended canonical}
}}
\put(5131,-1006){\makebox(0,0)[lb]{unknown territory}%
}
\put(4726,-6406){\makebox(0,0)[lb]{unknown territory}%
}
\put(2656,-2626){\rotatebox{270.0}{\makebox(0,0)[lb]{unknown territory}%
}}
\end{picture}%
\end{center}
We have pointed out that (derived) accessible algebras are especially suitable for a spectral analysis. Still the situation is not perfect since accessible algebras may be isospectral without being derived equivalent. This poses the question to isolate interesting subclasses $\Cc$ like path algebras of Dynkin quivers, or path algebras of stars, or canonical algebras with at most three weights where $\Cc$ has the \emph{separation property}. This means that two members of $\Cc$ with the same spectrum must be derived equivalent. It is an open question whether the class $\Cc$ of extended canonical algebras with three weights has this separation property. We conjecture that this is indeed the case, and moreover that the weight type can always be recovered from the Coxeter polynomial of an extended canonical algebra, regardless what the number of weights is.

Of course, all such questions are special instances of the general question to understand which properties of the representation theory of a finite dimensional algebra can be recovered from its spectral properties.

\end{document}